# Algebraic generality *vs* arithmetic generality in the controversy between C. Jordan and L. Kronecker (1874).

**Frédéric Brechenmacher (1).**

## *Introduction.*

Throughout the whole year of 1874, Camille Jordan and Leopold Kronecker were quarrelling over two theorems. On the one hand, Jordan had stated in his 1870 *Traité des substitutions et des équations algébriques* a canonical form theorem for substitutions of linear groups; on the other hand, Karl Weierstrass had introduced in 1868 the elementary divisors of non singular pairs of bilinear forms ($P$,$Q$) in stating a complete set of polynomial invariants computed from the determinant $|P+sQ|$ as a key theorem of the theory of bilinear and quadratic forms. Although they would be considered equivalent as regard to modern mathematics (2), not only had these two theorems been stated independently and for different purposes, they had also been lying within the distinct frameworks of two theories until some connections came to light in 1872-1873, breeding the 1874 quarrel and hence revealing an opposition over two practices relating to distinctive cultural features. As we will be focusing in this paper on how the 1874 controversy about Jordan's algebraic practice of canonical reduction and Kronecker's arithmetic practice of invariant computation sheds some light on *two conflicting perspectives on polynomial generality* we shall appeal to former publications which have already been dealing with some of the cultural issues highlighted by this controversy such as tacit knowledge, local ways of thinking, internal philosophies and disciplinary ideals peculiar to individuals or communities [Brechenmacher 200?a] as well as the different perceptions expressed by the two opponents of a long term history involving authors such as Joseph-Louis Lagrange, Augustin Cauchy and Charles Hermite [Brechenmacher 200?b].

---

1 F. BRECHENMACHER. Laboratoire de Mathématiques Lens (LML, EA2462). Fédération de Recherche Mathématique du Nord-Pas-de-Calais (CNRS, FR 2956). Université d'Artois (IUFM du Nord Pas de Calais). Faculté des Sciences Jean Perrin, rue Jean Souvraz S.P. 18, 62 307 Lens Cedex France.
Email: frederic.brechenmacher@euler.univ-artois.fr.

2 From the standpoint of modern algebra, Jordan's canonical form theorem for matrices with coefficients belonging to an algebraically closed field is equivalent to the elementary divisors theorem. An elementary divisor $(\lambda\text{-}a)^k$ corresponds to a $k$ by $k$ Jordan block, see [Gantmacher 1959]. See also the notes n° 48 and n°52 about Kronecker's and Frobenius' developments on invariant factores. We give below some examples of relations between three Jordan's canonical forms and three polynomial decompositions of the characteristic polynomial $|A\text{-}\lambda I|=(\lambda\text{-}1)^2(\lambda\text{-}2)^3(\lambda\text{-}3)$.

| | | | |
|---|---|---|---|
| Jordan's canonical forms. | $\begin{pmatrix}1&&&&&\\&1&&&&\\&&2&&&\\&&&2&&\\&&&&2&\\&&&&&3\end{pmatrix}$ | $\begin{pmatrix}1&&&&&\\&1&&&&\\&&2&1&&\\&&&2&1&\\&&&&2&\\&&&&&3\end{pmatrix}$ | $\begin{pmatrix}1&&&&&\\&1&&&&\\&&2&0&&\\&&&2&1&\\&&&&2&\\&&&&&3\end{pmatrix}$ |
| Elementary divisors. | *(λ-1), (λ-1), (λ-2), (λ-2), (λ-2), (λ-3)* | *(λ-1), (λ-1),* $(\lambda\text{-}2)^3$, *(λ-3)* | *(λ-1), (λ-1), (λ-2),* $(\lambda\text{-}2)^2$, *(λ-3)* |

In this purpose we shall frequently refer to a series of essays the historian Thomas Hawkins had devoted to the history of the theory of matrices in the 19th century and in which he argued that although "historians writing on this subject have tended to emphasize the role Arthur Cayley ([3]) […] there is much more to the theory of matrices – and to its history – than the formal aspect, i.e. the symbolical algebra of matrices. There is also a content: the concepts and theorems that make it a bona fide theory. An important part of that content is what I have termed spectral theory: the concept of an eigenvalue, the classification of matrices into types (symmetric, orthogonal, Hermitian, unitary, etc.), the theorems on the nature of the eigenvalues of the various types and, above all, those on the canonical (or normal) forms for matrices" [Hawkins 1975, p. 1]. In 1975, Hawkins related Cauchy's 1829 important memoir on the classification of conics and quadrics to the mechanical investigations of Lagrange and Pierre-Simon Laplace thereby considering what he referred to as "eigenvalue problems" throughout the 18th an 19th centuries ([4]). A few years later he presented what he considered as the "keystone" of this history: "the theory of canonical matrix forms" and laid special emphasis on Weierstrass' theorem of elementary divisors because "in [1868] and a preliminary memoir [1858] Weierstrass demonstrated more than theorems. He also demonstrated the possibility and desirability of a more rigorous approach to algebraic analysis that did not rest content with the prevailing tendency to reason vaguely in terms of the "general case". Weierstrass' memoir thus served as a paradigm for further research by his colleagues and students such as Kronecker and Frobenius" so that, according to Hawkins, "through his theory of elementary divisors [1868] and his influence upon other mathematicians, he [Weierstrass] more than anyone, was responsible for the emergence of the theory of matrices as a coherent, substantial branch of twentieth-century mathematics" [Hawkins 1977, p. 119]. Although Hawkins' work laid foundations for the investigations presented in this paper, we shall also develop new considerations in wondering about Kronecker and Jordan different perspectives on generality without focusing on the issues about the origins of abstract notions, theories, ways of reasoning and, more generally, structures, most authors have been dealing with while studying the history of linear algebra. The theoretical identity of the 1920-1930's "modern algebra" has often served as a lens for looking into the past, selecting relevant texts and authors, thereby giving structure to its own history while other identities that did not fit in this retrospective theoretical glance have stayed out of sight. The question therefore arises as to the *identities* and *significations* taken on by the *algebraic practices* that had been developed within various disciplines - such as mechanics, arithmetic or geometry - and had passed from one theory on to another before the time of unifying *algebraic theories* such as the 1930's theory of matrices when these *practices* would be seen as *methods* for the classification of similar matrices ([5]). As we will be looking into the

[3] In 1858 A. Cayley published *A memoir on the theory of matrices* in which he developed some symbolic laws of addition and multiplication on the "matrices" that had been introduced by Sylvester in 1851. See [Hawkins 1977] and [Brechenmacher 2006c] as well as M.J. Durand-Richard's paper n this volume devoted to some perspectives on generality in Cayley's mathematics.

[4] From the standpoint of modern algebra, Cauchy's 1829 memoir provided the first "general" proof that the eigenvalues of a symmetric matrix are real and that the corresponding quadratic form can be transformed into a sum of square terms (i.e. diagonalized) by means of an orthogonal transformation.

[5] This paper originated in a doctoral research devoted to the Jordan canonical form theorem (1870-1930) under the supervision of J. Dhombres [Brechenmacher 2006a]. Major sources of inspirations have been C. Goldstein'investigation s on the identities taken on by a theorem of Fermat [Goldstein 1995], H. Sinaceur's history of Sturm theorem [Sinaceur 1991], the works of C. Gilain on the fundamental theorem of algebra [Gilain 1991], the researches of G. Cifoletti on Pelletier and Gosselin algebraic practices [Cifoletti 1992] as well as the

1874 controversy our purpose will thus also be to discuss how Jordan's and Kronecker's internal philosophies of generality were related to some complex issues about the algebraic or non algebraic natures of the practices they were quarrelling about.

## *1. A polemical "general" theory.*

The controversy started in the winter 1873-1874 when two papers were successively read to the Academies of Paris and Berlin and, although it would stay at the beginning on the private level of a correspondence between Jordan and Kronecker, it would turn in spring and summer into a public quarrel and to the publication of a train of notes and papers ([6]). The controversy was originally caused by Jordan's ambition to reorganise the theory of bilinear forms through what he designated as the "simple" notion of "canonical form". Jordan's December 1873 note was actually the first contribution to the theory of bilinear forms to be published out of Berlin and the issue at stake was the organisation of what used to be a local field of research limited to a few Berliners now in the process of turning into an international theory. In 1866, two papers published by Elwin Christoffel and Kronecker in Crelle's Journal had laid the foundations of a theory whose main problem was the characterisation of bilinear forms – given two forms *P* and *P*', find the necessary and sufficient conditions under which *P* can be transformed into *P'* by using linear substitutions – and whose methods were to look for invariants computed from the forms coefficients and which would be unaltered by linear transformations. It was actually the problem of the simultaneous transformations of two forms *P* and *Q* which would shortly become the main question of the theory. Although the determinant of the "network" ("Schaaren") $P+sQ$ was a polynomial invariant, the roots of the characteristic equation $|P+sQ|=0$ would provide a complete set of invariants only under the condition that no multiple roots existed ([7]). The general resolution of this problem, which had originally been tackled by Kronecker in 1866 in the context of his researches on abelian functions, regardless of whether $|P+sQ|$ factored into distinct linear factors or not, had been given in two successive memoirs published in 1868. While Weierstrass introduced a complete set of invariants computed from a comparison of the algebraic decompositions of the determinant $|P+sQ|$ and its successive minors, therefore giving a resolution to the non singular case where $|P+sQ|$ did not vanish identically, Kronecker's forthcoming paper was devoted to the singular case. Since then, Weierstrass' elementary divisors theorem had become the main result of the theory of bilinear forms. Jordan thus stroke at the heart of a theory when he claimed in December 1873 that "the problem of the simultaneous reduction of two functions *P* and *Q* is identical to the problem of the reduction of a linear substitution to its canonical form" [Jordan 1873, p. 1487] and to the theorem he had stated in his 1870 *Traité des substitutions et des equations algébriques*. Jordan proposed a reorganisation in three questions of the theory by appealing to a "general method" consisting in "reducing" "bilinear polynomials" to some "canonical forms" relating to the types of substitutions group considered.

> It is known that any bilinear polynomial

---

long term perspectives given on some fluctuations of mathematical elaborations and on the history of fractions in the collective book coordinated by P. Benoit, K. Chemla et J. Ritter [Benoit and al. 1992].

[6] For a full study of the 1874 controversy and of the unpublished correspondence between Jordan and Kronecker see [Brechenmacher 200?a].

[7] For instance, the bilinear forms $B(X,U)=ux+uy+vy$ and the identity $I=ux+vy$ both have a single eigenvalue equal to 1 but these two forms are not equivalent one two another.

$$P= \Sigma A_{\alpha\beta}x_\alpha y_\beta \ (\alpha=1,2,...,n, \beta=1,2,...,n)$$

can be reduced to its canonical form $x_1y_1+...+x_my_m$, by linear transformations applied to the two sets of variables $x_1,..,x_n$, $y_1,..,y_n$. We now consider the following questions among the many of this kind :

1. To reduce a bilinear polynomial $P$ to its canonical form by applying orthogonal substitutions on the two systems of variables $x_1,...,x_n$ ; $y_1,...,y_n$.
2. To reduce a bilinear polynomial $P$ to its canonical form by using some substitutions simultaneously on the $x$'s and the $y$'s.
3. To reduce simultaneously two bilinear polynomials $P$ and $Q$ to a canonical form by using some linear substitutions on each set of variables individually. [Jordan 1873, p. 7-11, translation F.B.] ([8]).

In the 1870's, not only did many recent applications herald the major role the theory of bilinear forms would play in the following decades, this theory was also giving a new "homogeneous" and "general" treatment to different problems referring to various theories developed throughout the 19th century ([9]). The above quoted note Jordan communicated to the Parisian Academy on the 22nd of December alluded to geometry and Cauchy's results on the principal axis of conics and quadrics (first question), to the arithmetic of quadratic forms relating especially to the works of Carl Gauss and Hermite (second question) as well as to analytical mechanic and the solution given by Lagrange to systems $PY''+QY=0$ of linear differential equations with constant coefficients (third question) ([10]). As we shall see in greater detail later, it was the *general* solutions they were both able to give to these problems which had been tackled in the past by Lagrange, Cauchy or Hermite that lead to a connection between the two theorems stated independently by Jordan and Weierstrass and which therefore prompted Jordan's 1873 irruption in the theory of bilinear forms.

---

[8] On sait qu'il existe une infinité de manières de ramener un polynôme bilinéaire $P= \Sigma A_{\alpha\beta}x_\alpha y_\beta$ $(\alpha=1,2,...,n, \beta=1,2,...,n)$ à la forme canonique $x_1y_1+...+x_my_m$, […] par des transformations linéaires opérées sur les deux systèmes de variables $x_1,..,x_n$, $y_1,..,y_n$. Parmi les diverses questions de ce genre que l'on peut se proposer, nous considérons les suivantes : 1. Ramener un polynôme bilinéaire $P$ à une forme canonique simple par des substitutions orthogonales opérées les unes sur $x_1,...,x_n$, les autres sur $y_1,...,y_n$. 2. Ramener $P$ à une forme canonique simple par des substitutions linéaires quelconques opérées simultanément sur les $x$ et les $y$. 3. Ramener simultanément à une forme canonique deux polynômes $P$ et $Q$ par des substitutions linéaires quelconques, opérées isolement sur chacune des deux séries de variables.

[9] In modern terms, bilinear forms played for a long time a role similar to the one matrices would be playing in XXth century linear algebra. Developed at the beginning as a local field of research limited to a few geometers from Berlin such as Christoffel, Kronecker and Weierstrass, the theory of bilinear forms expanded to an international theory between 1874 and 1880 thanks to the numerous applications developed in geometry (Klein 1868), the theory of quadratic forms (Kronecker 1874, Darboux 1874), various problems related to systems of differential equations (Jordan 1871-1872) as Fuchs equations (Hamburger 1872, Jordan 1874) or Pfaff's problem (Frobenius and Darboux 1875-1880). See [Hawkins 1977] and [Brechenmacher 2006a]. On Christoffel's works see especially [Mawhin 1981].

[10] From the standpoint of modern algebra, the canonical form $x_1y_1+...+x_my_m$ characterizes the equivalence classes of square matrices for the equivalence relation $(ARB \Leftrightarrow \exists P,Q \in GL_n(\hat{E}), PAQ=B)$. The first problem concerns the similarity relation for orthogonal matrices $(ARB \Leftrightarrow \exists P \in O(\hat{E}), P^{-1}AP=B)$. The second problem relates to the congruence relation for square matrices $(ARB \Leftrightarrow \exists P \in GL_n(\hat{E}), {}^tPAP=B))$. The third problem focuses on the equivalence of pairs of matrices $(A, B)$ which is of importance for systems of differential equations with constant coefficients $AY''+BY=0$. In the particular case where $B=I$, the equivalence relation for pairs $(A,I)$ corresponds to the similarity relation $B=P^{-1}AP$. In the above quotation, Kronecker pointed out that the first and second problems might be deduced from the third one because they might be respectively considered as relating to the congruence of pairs $(A,I)$ and to the equivalence of pairs $(A, {}^tA)$. For some mathematical complements on the equivalence of pairs of matrices see [Dieudonné 1946].

> […] the third [problem has already been handled] by M. Weierstrass, the solutions given by the geometers from Berlin are nevertheless incomplete as some exceptional cases have been left aside despite their significance. Their analysis is also quite difficult to follow – especially M. Weieerstrass' one – and we shall therefore suggest an extremely simple new method that holds no exception […]. We will show that the problem of the simultaneous reduction of two functions *P* and *Q* is identical to the problem of the reduction of a linear substitution to its canonical form. [*Ibidem*] ([11]).

In the paper he communicated to the Berlin Academy on the 19th of January, while responding to Jordan's claims for the greater « simplicity » and « generality» of his methods of canonical reduction as opposed to Weierstrass' invariant computations, Kronecker did not only reject the originality and validity of these methods but also the theoretical organisation relating to them. Recalling that as soon as 1868 Weierstrass and he had organised the theory around the sole problem of the characterisation of pairs of forms, he questioned the relevance of Jordan's distinction between three problems of canonical reductions.

> In M. *Jordan*'s paper, "On bilinear forms", the solution to the first problem is not original, the solution to the second problem is incorrect and the solution to the third one is incomplete. In addition to that, the third problem, which actually includes the two others as particular cases, has been completely solved by M. *Weierstrass*' work of 1868 and by my additional contribution to it. Unless I am very much mistaken, there are serious grounds for questioning M. *Jordan*'s priority in the invention of his results, should they even be correct. [Kronecker 1874b, p. 1181, translation F.B.] ([12]).

During winter, Kronecker would develop his views on the organisation of the theory of bilinear forms in communicating monthly papers to the Academy of Berlin. Meanwhile, he and Jordan engaged a private correspondence with the aim of settling the quarrel of priority bred by the connections which had arisen between the canonical form and the elementary divisors theorems ([13]). Despite its function of scientific communication and although it would lead Jordan to recognise a partial anteriority of Kronecker and Weierstrass on some of his results as well as to grasp some tacit knowledge peculiar to the Berliners' practices

---

[11] […] le troisième [problème a déjà été traité] par M. Weierstrass ; mais les solutions données par les éminents géomètres de Berlin sont incomplètes, en ce qu'ils ont laissé de côté certains cas exceptionnels qui, pourtant, ne manquent pas d'intérêt. Leur analyse est en outre assez difficile à suivre, surtout celle de M. Weierstrass. Les méthodes nouvelles que nous proposons sont, au contraire, extrêmement simples et ne comportent aucune exception […]. La réduction simultanée de deux fonctions *P* et *Q* est un problème identique à celui de la réduction d'une substitution linéaire à sa forme canonique

[12] Dans le Mémoire de M. *Jordan [...],* la solution du premier problème n'est pas véritablement nouvelle ; la solution du deuxième est manquée, et celle du troisième n'est pas suffisamment établie. Ajoutons qu'en réalité ce troisième problème embrasse les deux autres comme cas particuliers, et que sa solution complète résulte du travail de M. *Weierstrass* de 1868 et se déduit aussi de mes additions à ce travail. Il y a donc, si je ne me trompe, de sérieux motifs pour contester à M. *Jordan* l'invention première de ses résultats, en tant qu'ils sont corrects [...].

[13] This correspondence which is kept in the archives of the Ecole Polytechnique has been edited in [Brechenmacher 2006a].

([14]), the correspondence would fail to reach agreement on a mathematical ground. On the one hand, Kronecker failed to bring Jordan round to his own ideas on the structure of the theory of bilinear forms, on the other hand Jordan did not succeed in convincing Kronecker of his "natural right" to claim the genuine originality of his method which, he said, "has rather brought to light than disparaged Weierstrass' result in showing the resolution it implicitly gave to a fundamental problem in linear substitutions theory which, to my opinion, is a much more fertile theory than the algebraic theory of forms of the second order" [Jordan to Kronecker, January 1874, translation F.B.]. In spring, the controversy would go public again and reach its climax with two successive communications of Jordan and Kronecker to the Academy of Paris. The quarrel of priority would then turn into an opposition over two theories (group theory vs theory of bilinear forms) and two disciplines (algebra vs arithmetic) as well as over two practices (canonical reduction vs invariant computation) relating to conflicting philosophies of generality.

## *2. The generality of Weierstrass' theorem as marking a rupture in the history of algebraic reasonings.*

As we shall see in this section, Kronecker associated Weierstrass' theorem with an ideal of generality which revealed itself on the occasion of his criticisms of the formal nature he imputed to Jordan's canonical reduction. Because he had designated as canonical forms three different algebraic expressions relating to his classification in three problems, Jordan was accused of resorting to a notion without any "general relevance" nor "objective content", therefore mixing up the "formal aspects" of some "means of action" (canonical forms) with the "true subject of investigation" and its "content". Although "normal forms" similar to Jordan's canonical form were being used by Kronecker himself; to the latter's opinion, resorting to such algebraic expressions was legitimate provided that they would remain in their proper places of "methods" as opposed to the "notions" relating to the "other disciplines" - such as arithmetic - it was algebra's duty to serve ([15]). On the contrary, mixing up algebraic methods with notions and theoretical issues would lead to mistake a mere "formal" development for a "general" and "uniform" presentation. Kronecker thus mocked Jordan's claims for the greater simplicity of his canonical reduction as a naïve simplism which contented itself with the illusionary generality and uniformity of a formal development.

---

[14] In the analysis of the correspondence we have published in [Brechenmacher 200?a], some tacit knowledge peculiar to the Berliners' network have been highlighted. Until he received some detailed explanations thanks to his correspondence with Kronecker, Jordan had some difficulties in understanding the implicit properties of determinants Weierstrass and Kronecker resorted to for computing invariants as well as the connections between the various papers published during the 1860's (especially because some paper were not explicitly mentioning "bilinear forms" in their titles and Jordan believed them to be only devoted to "quadratic forms"). For instance, when he published his first note in December 1873, Jordan was unaware of Christoffel and Kronecker works on the bilinear forms and he did not realize the relation between Kronecker's 1868 memoir on "quadratic forms" and Weierstrass' 1868 paper despite the fact that these two papers had been conceived as the "two parts" of a same development and had thus been published successively in *Crelle's* Journal.

[15] For instance, in order to demonstrate that two non singular pairs of bilinear forms could be transformed one into the other, Weierstrass proved in 1868 that both could be linearly transformed to what Kronecker would refer to as a normal form similar to Jordan's canonical form. But neither Kronecker nor Weierstrass would state a theorem about such normal forms which were not the purpose of their investigations. On Weierstrass' 1868 demonstration see [Hawkins 1977] and [Brechenmacher 2006a].

Should such general expressions be found, simplicity and generality motivations might eventually lead us to designate them under a common name of "canonical forms". But shall we refuse to content ourselves with the mere formal aspects the most recent algebraic works have so often been pushing forward – with the most unlikely benefits for the cause of science -, we shall not neglect to justify by internal grounds the relevance of such canonical forms. These so called canonical or normal forms are actually established only because of the orientation given to the investigation and we thus shall consider them only as means of action as opposed to the aims of the research [….]. It is not surprising that while looking for a completely general and uniform exposition the author [Jordan] had been forced to introduce new principles in the paper already referred to [Jordan 1874a]. We should have been very surprised if, on the opposite and in accordance to Jordan's claims ("The new methods we are proposing are, on the contrary, extremely simple…" "It is plain to see that we can transform…."), the simplest means would have been sufficient ([16]). [Kronecker 1874a, p.367, translation F.B.].

Kronecker then developed his views on the opposition of two kinds of "generality". Taking as a starting point the "fault" of using as a denominator an algebraic expression that might vanish he picked out Jordan's paper [1874a] ([17]), he blamed the "so called generality" of formal expressions which would actually lose any meaning in some singular cases.

We are getting used to discover some new difficulties – especially in algebraic issues - as soon as we free ourselves from the restriction to such cases we would be tempted to consider as the general ones following our usual custom. As soon as we force our way through the surface of this so called generality which actually excludes any particularity, we penetrate inside the true generality covered with all its singularities and it is usually not until then that we come to face the true difficulties of a subject as well as the many new phenomenons and point of views concealed at its depth [*Ibidem*] ([18]).

---

[16] Nachträglich, wenn dergleichen allgemeine Ausdrücke gefunden sind, dürfte die Bezeichnung derselben als canonische Formen allenfalls durch ihre Allgemeinheit und Einfachheit motiviert werden können ; aber wenn man nicht bei den bloss formalen Gesichtspunkten stehen bleiben will, welche –gewiss nicht zum Vortheil der wahren Erkenntnis- in der neueren Algebra vielfach in den Vordergrund getreten sind, so darf man nicht unterlassen, die Berechtigung der aufgestellten canonischen Formen aus inneren Gründen herzuleiten. In Wahrheit sind überhaupt die so genannten canonischen oder Normalformen lediglich durch die Tendenz der Untersuchung bestimmt und daher nur als Mittel, nicht aber als Zweck der Forschung anzusehen. […] Dass sich aber für eine zugleich einheitliche und ganz allgemeine Entwickelung, wie sie in der ben erwähten Arbeit gegeben ist, gewisse neue Principien als nöthig erwiesen, kann durchaus nicht befremden, und es wäre im Gegentheil zu verwundern, wenn wirklich den *Jordan*'schen Behauptungen gemäss ("Les méthodes nouvelles que nous proposons sont, au contraire extrêmement simples..." "On voit par une discussion très simple, que l'on peut transformer...") die allereinfachsten Mittel dazu ausreichen sollten.

[17] Jordan promptly corrected this mistake which was of no consequence for his theoretical organisation, see [Brechenmacher 2006a, p. 689].

[18] Denn man ist es gewohnt –zumal in algebraischen Fragen- wesentlich neue Schwierigkeiten anzutreffen, wenne man sich von der Beschränkung auf diejenigen Fälle losmachen will, welche man als die allgemeinen zu bezeichnen pflegt. Sobald man von der Oberfläche der sogenannten, jede Besonderheit ausschliessenden Allgemeinheit in das Innere der wahren Allgemeinheit eindringt, welche alle Singularitäten mit umfasst, findet man in der Regel erst die eigentlichen Schwierigkeiten der Untersuchung, zugleich aber auch die Fülle neuer Gesichtpunkte und Erscheinungen, welche sie in ihren Tiefen enthält.

To Kronecker's mind, Weierstrass' theorem was a model of "true generality", the general and homogeneous resolution the elementary divisors were giving to the problem of the characterisation of bilinear forms contrasted with the "inadequate results" of the "so called general" methods that had been sporadically developed "for over a century" with little attention given to difficulties that might be caused by singularities such as equal factors occurring in the polynomial determinant $S=|P+sQ|$. Kronecker was implicitly making reference to the works of authors such as Lagrange, Sturm, Cauchy or Carl Gustav Jacobi whose algebraic practices were blamed for their tendency to think in terms of the "general" case where $S=0$ had no multiple roots. As we shall see in greater detail later, such practices had been elaborated in investigating the symmetric case of "networks of quadratic forms", they resorted to "general" polynomial expressions relating to polynomial factorisations of $S$ and its sub determinants $P_{1i}$ (*i.e.* what we would now call matrix minors lying on the intersection of the first row and $i^{th}$ column subsets). Assigning specific values to the symbols involved in the general expression (*) $\frac{P_{1i}}{\frac{S}{x-s_j}}(x)$ (which gave the systems $(x_i^{s_j})_{1\leq i\leq n}$ of solutions related to the root $s_j$ : $x_i^{s_j}=\frac{P_{1i}}{\frac{S}{x-s_j}}(s_j)$ ) ([19]), would hence cause

[19] Given a symmetric matrix $A$, the coordinates of its eigenvectors are given by not equal to zero columns of the adjoint matrix of the characteristic matrix $A$-$xI$ (*i.e.* the matrix of cofactors). For instance, the matrix :

$$A=\begin{pmatrix}1 & -1 & 0\\ -1 & 2 & 1\\ 0 & 1 & 1\end{pmatrix}$$

The characteristic equation and minors are given by : $S=-x(3-x)(1-x)$, $P_{11}(x)=(1-x)(2-x)-1$, $P_{12}(x)=(1-x)$ et $P_{13}=-1$. The expression $((1-x)(2-x)-1,\ 1-x,\ -1)$ gives the polynomial coordinates of any eigenvector (while $|S|$ gives the square of its norm), *e.g.* for the eigenvalue $s_1=1$, $\frac{S}{x-1}=x(3-x)$,

$$x_1^{s_1}=\frac{P_{11}}{\frac{S}{x-1}}(1)=-\frac{1}{2},\ x_2^{s_1}=\frac{P_{12}}{\frac{S}{x-1}}(1)=0,\ x_3^{s_1}=\frac{P_{11}}{\frac{S}{x-1}}(1)=-\frac{1}{2},$$

The coordinates of the normed eigenvector relating to the eigenvalue 1 are $(1/\sqrt{2},\ 0,\ 1/\sqrt{2})$.
Similarly for the eigenvalues $s_2=0$,

$$x_1^{s_2}=\frac{1}{3},x_2^{s_2}=\frac{1}{3},x_3^{s_2}=-\frac{1}{3}$$

and $s_3=3$,

$$x_1^{s_3}=\frac{1}{6},x_2^{s_3}=\frac{2}{6},x_3^{s_3}=-\frac{1}{6}$$

The corresponding normed vectors are $(1/\sqrt{3},\ 1/\sqrt{3},\ -1/\sqrt{3})$ et $(1/\sqrt{6},\ 2/\sqrt{6},\ -1/\sqrt{6})$.
Interpretation for quadratic forms : la form associated to $A$ in the canonical basis of Ë$^3$,

$$A(x_1,x_2,x_3)=x_1^2-2x_1x_2+x_2^2+2x_2x_3+x_3^2=1.X_1^2+0.X_2^2+3.X_3^2$$

with $(X_1,X_2,X_3)=(x_i^{s_j})_{i,j=1,2,3}A(x_i^{s_j})_{i,j=1,2,3}^{-1}$ where $(x_i^{s_j})_{i,j=1,2,3}$ is the change of orhonormal basis matrix.

$$(x_i^{s_j})_{i,j=1,2,3}=\begin{pmatrix}\frac{1}{\sqrt{2}} & \frac{1}{\sqrt{3}} & \frac{1}{\sqrt{6}}\\ 0 & \frac{1}{\sqrt{3}} & \frac{2}{\sqrt{6}}\\ \frac{1}{\sqrt{2}} & -\frac{1}{\sqrt{3}} & -\frac{1}{\sqrt{6}}\end{pmatrix}$$

Interpretation for (symmetric) bilinear forms :

$$A(x,y)=x_1y_1-x_1y_2-y_1x_2+x_2y_2+x_2y_3+y_3x_2+x_3y_3=1.X_1Y_1+0.X_2Y_2+3.X_3Y_3.$$

difficulties as it would lead to $\frac{0}{0}$ expressions if common roots should occur between $S$ and its successive sub determinants. In 1858, Weierstrass had proven that, in the quadratic case, each root of multiplicity $p$ of $S(x)=0$ had to be a root of $P_{1i}(x)=0$ of a multiplicity at least equal to $p-1$, the expression (*) were therefore defined regardless of the multiplicities of the roots ([20]). To Kronecker's mind the idea of elementary divisor took root in 1858 when Weierstrass had shifted the emphasis from roots multiplicity toward the investigation of the links between the polynomial decompositions of the successive determinant extracted from $S$ thereby generalising the inertia law of single quadratic forms in stating that it was possible to transform simultaneously two quadratic forms ($A=\sum_{i=1}^{n} A_{ij}x_i x_j$ $B=\sum_{i=1}^{n} B_{ij}x_i x_j$)- $A$ being definite positive - into sums of square terms ($\sum_{i=1}^{n} X_i^2$ , $\sum_{i=1}^{n} s_i X_i^2$ ) by means of what we would designate today as orthogonal substitutions. Because it proved that the roots $s_1, s_2, ..., s_n$ of $S=0$ were real whatever their multiplicity, the "true generality" of Weirstrass' 1868 statement heralded the 1860's generalisations to bilinear forms ([21]).

> [Such a situation] shows up in each of the few algebraic questions which have been tackled in minute details such as the theory of networks of quadratic forms whose main characteristics have been developed in this paper. Because for a long time no one had been daring to venture beyond the condition that only unequal factors would occur in the determinant [$|P+sQ|$], only inadequate results had been achieved in the classical problem of simultaneous transformations of two quadratic forms and the many sporadic attempts to deal with this subject in the past century had been ignoring its true nature. Freeing itself from this condition, Weierstrass' 1858 work had resulted in an elevated conception in providing a complete resolution of the case when only simple elementary divisors occurred. […] Despite this first step, it was not until the publication of Weierstrass' 1868 paper that the fully general notion of elementary divisor was introduced thereby shedding a new light on the theory of networks of forms whatever the situation with the sole restriction that the determinant should not vanish. After I had gone through this last restriction in developing the notion of elementary divisor of general elementary networks, the brightest light began to shine on numbers of new algebraic forms and, through this complete resolution of the

[20] As opposed to the general case of bilinear forms, in the quadratic case Weierstrass had considered in 1858 the matrices were symmetrical and therefore diagonalizable.

[21] According to the inertia law as it had been introduced by Hermite and Sylvester in the 1850's, it was possible to transform a quadratic form $A$ into a sum of squares $A=\Delta_{n-1}X_1^2+\frac{\Delta_{n-2}}{\Delta_{n-1}}X_2^2+...+\frac{\Delta}{\Delta_1}X_n^2$, where $\Delta, \Delta_1, \Delta_2, ..., \Delta_{n-1}, 1$ were the principal minors of $A$. In 1874 G. Darboux explicitly generalised the inertia law to polynomial quadratic forms $A+sI$ [Darboux 1874, p. 367]. Darboux's paper was published next to Jordan's controversial one [1874a] in *Liouville's Journal*. Both proposed new demonstrations on Weierstrass' theorems: while Jordan focused on the 1868 result on pairs of bilinear forms, Darboux worked out new perspectives on the 1858 theorem on quadratic forms in developing the methods introduced by Hermite and Sylvester in the 1850's. Darboux was merely concern with the geometry of surfaces and his method was inserted by Gundelfinger in the third edition of Hesse's analytical geometry. For a detailed description of this method see [Drach and Meyer 1907].

> question, higher ideas had been gained on a theory of invariants developed in its true generality [*Ibid.*] ([22]).

Combining mathematical and historical arguments, Kronecker was the first to stress a history of what the historian T. Hawkins would refer to in the 1970's as "generic reasoning" in the 16th-19th centuries algebra.

> The generality of the method of analysis had been viewed as its great virtue since its inception. Thus Viète stressed that the new method of analysis "does not employ its logic on numbers – which was the tediousness of the ancient analysts – but uses its logic through a logistic which in a new way has to do with species" […]. Analysis became a method for reasoning with, manipulating, expressions involving symbols with "general" values and a tendency developed to think almost exclusively in terms of the "general" case with little, if any, attention given to potential difficulties or inaccuracies that might be caused by assigning certain specific values to the symbols. Such reasoning with "general" expressions I shall refer to for the sake of brevity as *generic reasoning.* [Hawkins 1977, p. 122].

The opposition between generic and general reasoning in the ways the characteristic equation's multiple roots were being handled structured a long term history of algebra from the mechanical systems studied by Lagrange and Laplace in the 18th century to the 19th century when some mathematicians such as Cauchy and Weierstrass concerned with raising the standards of rigor came to reject the legitimacy of generic reasoning:

> […] neither of them [Lagrange and Laplace] had pursued the study of the solutions of systems of linear differential equations with sufficient care to justify their claim [that the characteristic roots $\lambda$ must be real]. They had no difficulty treating such a system when the characteristic roots are distinct, but their analysis of the case of multiple roots was inadequate. Given the generic tendency of their analytical methods, it is noteworthy that they considered the case at all. […] Weierstrass recognition of the questionable nature of their claims formed the starting point of the investigation that culminated in his theory of elementary divisors. […] In [1868] and a preliminary memoir [1858] Weierstrass demonstrated more than theorems. He also demonstrated the possibility and desirability of a

[22] Diess bewährt sich durchweg in den wenigen algebraischen Fragen, welche bis in alle ihre Einzelheiten vollständig durchgeführt sind, namentlich aber in der Theorie der Schaaren von quadratischen Formen, die obein in ihren Hauptzügen entwickelt worden ist. Denn so lange man es nicht wagte, die Voraussetzung fallen zu lassen, dass die Determinante nur ungleiche Factoren enthalte, gelangte man bei jener bekannten Frage der gleichzeitigen Transformation von zwei quadratischen Formen, welche seit einem Jahrhundert so vielfach, wenn auch meist blos gelegentlich, behandelt worden ist, nur zu höchst dürftigen Resultaten, und die wahren Gesichtpunkte der Untersuchung blieben gänzlich unerkannt. Mit dem Aufgeben jener Voraussetzung führte die *Weierstrass*'sche Arbeit vom Jahre 1858 schon zu einer höheren Einsicht und namentlich zu einer vollständigen Erledigung des Falles, in welchem nur einfache Elementartheiler vorhanden sind. Aber die allgemeine Einführung dieses Begriffes der Elementartheiler, zu welcher dort nur ein vorläufiger Schritt gethan war, erfolgte erst in der *Weierstrass*'schen Abhandlung vom Jahre 1868, und es kam damit ganz neues Licht in die Theorie der Schaaren für den Fall beliebiger, doch von Null verschiedener Determinanten. Als ich darauf auch diese letzte Beschränkung abstreifte und aus jenem Begriffe der Elementartheiler den allgemeineren der elementaren Schaaren entwickelte, verbreitete sich die vollste Klarheit über die Fülle der neu auftretenden algebraischen Gebilde, und bei dieser vollständigen Behandlung des Gegenstandes wurden zugleich die wertvollsten Einblicke in die Theorie der höheren, in ihrer wahren Allgemeinheit aufzufassenden Invarianten gewonnen.

> more rigorous approach to algebraic analysis that did not rest content with the prevailing tendency to reason vaguely in terms of the "general" case. [Hawkins 1977, p. 123].

To both Kronecker's and Hawkins' views, Weierstrass' theorem marked a rupture in the history of algebra. Not only did it achieve a rigorous development as opposed to the generic nature of the reasonings of the past, it also resulted in a homogeneous resolution which proved that it was not necessary to resort to particular non algebraic arguments in order to handle the singular cases which restricted the range of validity of general algebraic expressions. In his work on the history of the spectral theory of matrices, Thomas Hawkins has displayed some evidences of such a rupture in the evolution of Kronecker's mathematical work. In a paper devoted to theta functions he published in 1866, Kronecker was still looking upon the occurrence of multiple roots as a "singular" case he had to tackle by appealing to arguments peculiar to the context of theta functions because his "general" algebraic method failed. In demonstrating the possibility of developing truly general and homogeneous methods, Weierstrass' invariants gave rise to new ideals of generality and not only would such ideals play a major role in the development of the theory of bilinear forms in Berlin, according to Hawkins they would also mark a major stage in the history of the 1930's "spectral theory of matrices".

> I would suggest that, insofar as anyone deserves the title of founder of the theory of matrices, it is Weierstrass. […] His theory of elementary divisors provided a theoretical core, a substantial foundation, upon which to build. His work demonstrated the possibility of dealing by the methods of analysis with the non-generic case, thereby opening up a whole new world to mathematical investigation, a world that his colleagues and students proceeded to explore. […] On motivational force common to the entire 19th century was a concern for a more rigorous level of reasoning in mathematics. […] A concern for higher standards of reasoning was a driving force behind Weierstrass' work and also behind that of Cauchy and Dirichlet which preceded it and behind that of Kronecker and Frobenius whoch succeeded it. The rise of the theory of matrices was directly related to the fall of the generic approach to algebraic analysis. A concern for rigor did not mark the end of the creative development of the theory but its beginning [Hawkins 1977, p. 157-159].

It was however in the very special context of a controversy that Kronecker came to emphasise the generic nature of some algebraic reasonings of the past century ([23]). The question therefore arises as to the different views that were being held by Jordan on what his opponent referred to as the "history of the theory of networks of quadratic forms". For the purpose of investigating how Jordan came to develop some connections between his researches on substitutions groups and the theory of forms we shall develop a complementary approach to the work of Thomas Hawkins on the generic nature of algebraic *reasoning* by focusing on the different significations taken on by an algebraic *practice* within various disciplinary frameworks –such as mechanics, geometry or arithmetics- before 1874 when this practice would be inserted in the two methods Jordan and Kronecker would rely to in their opposition on the organisation of a "general" theory of bilinear forms.

[23] On the construction of history by mathematical texts see [Goldstein 1995], [Dhombres 1998], [Cifoletti 1992 and 1995] and [Brechenmacher 2006c].

### 3. *A mechanical misinterpretation related to an algebraic practice dating back to the time of Lagrange.*

In a note he had addressed in 1870 to the geometers of the Academy of Paris, the astronomer Antoine Yvon-Villarceau pointed out a mistake in the classical method "for integrating the equations of a rotating solid body under the action of gravity" which had been "introduced by the illustrious author of the *Mécanique analytique* for the special case of the small oscillations of a loaded string whose equilibrium is slightly disturbed while one its end remains in position.". In 1766 Lagrange had devised a "general method" for the "general case" involving an arbitrary (finite) number of masses – as opposed to the particular case of a string loaded by two or three masses that had already been tackled by Jean d'Alembert - and hence for the resolution of a system of *n* linear differential equations with constant coefficients.

$$(a)\begin{cases} \frac{d^2y'}{dt^2}+A'y'+B'y''+C'y'''+...+N'y^{(n)}=0 \\ \frac{d^2y''}{dt^2}+A''y'+B''y''+C''y'''+...+N''y^{(n)}=0 \\ \frac{d^2y'''}{dt^2}+A'''y'+B'''y''+C'''y'''+...+N'''y^{(n)}=0 \\ \quad ... \\ \frac{d^2y^{(n)}}{dt^2}+A^{(n)}y'+B^{(n)}y''+C^{(n)}y'''+...+N^{(n)}y^{(n)}=0 \end{cases}$$

Making use of his principle of reduction of order, Lagrange had been looking for *n* particular solutions $y=E\sin(t\sqrt{K}+\varepsilon)$ related to the *n* independent equations $\frac{d^2y}{dt^2}+Ky=0$ whose combination would give the general solution. The *K*'s emerged as the roots of an algebraic equation of the $n^{th}$ order obtained through the "elimination" of the system's linear equations. As Villarceau illustrated it with a system of two equations, the problem had thus been reduced to an equation which turned out to characterise the mechanical system ([24]):

$$(a)\begin{cases} g\frac{d^2u}{dt^2}+a\frac{d^2s}{dt^2}+cu=0, \\ f\frac{d^2s}{dt^2}+a\frac{d^2u}{dt^2}+cs=0, \end{cases}$$

> [These methods] give the characteristic equation:
>
> $$\frac{c^2}{\rho^4}-(f+g)\frac{c}{\rho^2}+fg-a^2=0$$
>
> […] designating the absolute values of its roots as $\rho$ and $\rho'$, we find the following expressions for *s* and *u* :
>
> $$(f)\begin{cases} s=\alpha\sin(\rho t+\beta)+\alpha'\sin(\rho' t+\beta') \\ u=\frac{a\rho^2}{c-g\rho^2}\alpha\sin(\rho t+\beta)+\frac{a\rho'^2}{c-g\rho'^2}\alpha'\sin(\rho' t+\beta') \end{cases}$$
>
> [Yvon-Villarceau 1870, p. 763, translation F.B.].

---

[24] *u* and *s* are functions of *t* while *g*, *f*, *a* are constant coefficients. Note that in Villarceau's system *(a)*, as opposed to Lagrange's 1766 system quoted above, the pair *(a,c)* of coefficients of *u* and *s* in the second row is the mirror image of the coefficients *(c,a)* in the first row. In [Brechenmacher, 200?b] we have shown that this property of mechanical systems originated in the specific practice Lagrange had devised in 1766 for the problems of small oscillations.

It was therefore a necessary condition that this equation had *n* single roots for reducing the system to the *n* independent equations the roots were related to ($\frac{d^2u}{dt^2}+\rho u=0$ and $\frac{d^2u}{dt^2}+\rho' u=0$ in Villarceau's example). Villarceau's 1870 note aimed at criticising a mechanical interpretation dating back to Lagrange and according to which the presupposition of mechanical stability (the oscillation had to remain small) assured that only single roots could occur because multiple roots would cause unbounded oscillations as the "time *t* would be coming out of the sinus" and solutions would be taking the form $s=t\sin(\rho t+\beta)$ ([25]).

> I claim that this condition is not necessary for the oscillations to remain small. […] A homogeneous solid body of revolution oscillating around a point of its axis gives a very simple case where equal roots occur in the characteristic equation. As it is plain to see, the solid's oscillations will remain small provided that the initial impulsion is not too strong and that the solid's centre of gravity is chosen below its centre of suspension and not to far from the vertical axis passing through this point. [*Ibidem*, p. 765] ([26]).

Even though he pointed out some serious « deficiencies » in the « general » resolution for problems of small oscillations in some similar way as Kronecker would blame in 1874 the "so called generality" of algebraic expressions, Yvon-Villarceau did not aim at criticizing a tendency of generic reasoning, his purpose was to question a practice which consisted in combining some mechanical interpretations with the algebraic nature of the roots of a specific equation. Although Villarceau's intervention had been stemming from mechanical concerns such as the application of Lagrange's method to the long term perturbations of the parameters determining the planetary orbits, it brought up a theoretical question to the attention of the Academy's geometers. As the occurrence of multiple roots was in no contradiction to mechanical stability and hence to the possibility of "reducing" a system of *n* equations to *n* single independent equations, the question therefore arose as to the characterisation of such system which could be resolved into separates equations. This question prompted the publication by Jordan of two notes in 1871 and 1872. In 1871 Jordan applied his method of reduction of a linear substitution to a canonical form for reducing a system of differential equations with constant coefficients

$$\frac{dx_1}{dt}=a_1x_1+\ldots+l_1x_n,\ldots,$$
$$(I)\qquad \frac{dx_2}{dt}=a_2x_1+\ldots+l_2x_n,\ldots,$$

[25] From the standpoint of modern algebra, the stability of a system depends upon whether its matrix is diagonalisable or not. The inequality of the system's eigenvalues is a sufficient but not necessary condition. The mechanical systems studied by Lagrange are diagonalisable because they are symmetric.

[26] Je dis qu'il n'est pas nécessaire que cette condition soit remplie, pour que les petites oscillations se maintiennent. […] Voici un cas très simple, auquel correspondent des racines égales de l'équation caractéristique : c'est celui d'un corps solide, homogène et de révolution, oscillant autour d'un point pris sur son axe de figure. Chacun comprendra sans recourir au calcul, que la petitesse des oscillations est assurée dans ce cas, si le centre de gravité est, à l'origine du mouvement, au-dessous du centre de suspension, à une petite distance de la verticale passant par ce point, et si le mouvement oscillatoire initial est suffisamment faible.

.......

$$\frac{dx_n}{dt} = a_n x_1 + \ldots + l_n x_n,$$

to a reduced form regardless of the multiplicity of the roots.

$$(6)\ \frac{dy_1}{dt} = \sigma y_1,\ \frac{dz_1}{dt} = \sigma z_1 + y,\ \frac{du_1}{dt} = \sigma u_1 + z_1, \ldots,\ \frac{dw_1}{dt} = \sigma w_1 + v_1$$

Such a reduced form could then be integrated directly to yield $w_1 = e^{\sigma t}\psi(t)$, $v_1 = e^{\sigma t}\psi'(t), \ldots,$ $y_1 = e^{\sigma t}\psi^{r-1}(t)$, where $\psi(t)$ was some polynomial of degree $r-1$ (where $r$ is the number of "variables" in the series $y_1, \ldots, w_1$ [Jordan 1871, p.787]). He also characterised the systems which could be resolved into $n$ separates equations $\frac{dy_i}{dt} = \sigma y_i$ each of which could be directly integrated and which solutions involved no algebraic factors $\psi(t)$ by the necessary and sufficient condition that each characteristic root $K$ of multiplicity $\mu$ had to be a root of each of the minors of order $\mu-1$. In 1872, Jordan proved that Villarceau's mechanical systems verified this property because of their quadratic nature; thereby he extended the initial mechanical question not only to the theory of linear substitutions but also to the theory of quadratic forms.

> It is plain to see that the question of the reduction of the system (1) to its canonical form is identical to the well known following problem: *To have the angles of the variables of two quadratic forms* T *and* U *disappear at the same time* [Jordan 1872, p. 320, translation F.B.] ([27]).

In proving that the multiplicity of roots was of no relevance to the subject of mechanical stability, Jordan reached the same conclusion that Weierstrass had already been giving in 1858 (symmetric case) and in 1868 ([28]). It was thus thanks to a hundred-year old mechanical problem that a first connection between Jordan's and Weierstrass' theorems came to light. This connection was pointed out in 1873 by Meyer Hamburger who mixed Jordan's canonical reduction together with Weierstrass invariant computation with the aim of giving a "general" resolution to Fuchsian equations in the case of multiple roots [Hamburger 1873, p. 113] ([29]), thereby calling Jordan's attention to the work of Weierstrass on the theory of bilinear forms and prompting the ensuing controversy when Jordan pointed out in 1873 that the transformation of forms could be seen as the

---

[27] Il est clair que la question de la réduction du système (1) à la forme canonique (7) est identique à ce problème connu : *Faire disparaître les angles des variables à la fois dans les deux formes quadratiques* T *et* U.

[28] Although Weierstrass had already stated such a condition for the symmetric case in 1858, he communicated the same result Jordan had stated in 1872 in a communication to the Berlin Academy in 1875 he presented as an application of his theorem of elementary divisors (making therefore no reference to Jordan). From the standpoint of Weierstrass' 1868 theorem, Jordan's condition is equivalent to the necessary and sufficient condition stating that two bilinear forms $P$ and $Q$ might be transformed simultaneously into sums of square terms if and only if the elementary divisors of $pP+qQ$ are linear. If $P$ and $Q$ are real quadratic forms such that $pP+qQ$ is definite for some $p$ and $q$, then the elementary divisors are linear and $P$ and $Q$ can be transformed simultaneously into sums of square terms.

[29] In 1866 Fuchs had investigated differential equations of the form

$$\frac{d^n y}{dz^n} + p_1 \frac{d^{n-1} y}{dz^{n-1}} + \ldots + p_{n-1} \frac{dy}{dz} + p_n = 0$$

with variable coefficients $p_i = p_i(z)$. In 1873, Hamburger pointed out that Jordan's canonical reduction could be used to "simplify" the form Fuchs had given to fundamental systems of solutions in case of multiple roots. See [Hawkins 1977, p.147].

composition of linear substitutions and proved that his canonical reduction could be used to derive Weierstrass' theorem ([30]).

### *4. On generality and the algebraic identity of a practice carried by the equation to the secular inequalities in the Planetary Theory.*

It was thus because of their capacity to give a general resolution to some problems that had been handled throughout the 18th and 19th centuries that some identities between Jordan's and Weierstrass' theorems had arisen between 1870 and 1873. Both theorems were shedding some new light on the past as they were making some results of authors such as Lagrange or Cauchy appear incomplete *i.e.* limited to what would be considered from now on as the special case in which only single roots occurred. The 1874 controversy might therefore be considered as opposing two different ends given to a common history and for the purpose of a deeper understanding of the role played by such a history in the quarrel, a bibliographic research has been carried out, starting from the authors and texts Jordan and Kronecker referred to and working out systematically the references that appeared successively in this paper chase. This methodology resulted in a network covering the period 1766-1874 a simplified representation of which is given in annex ([31]). As it is plain to see by looking at the main knots appearing in the entanglement of bibliographic references and which point to the mechanical work of Lagrange as well as to Cauchy's analytical geometry, this network can neither be identified to a theory nor to the resolution of a single problem (even though there was a problem at the origin of the network, this problem was being considered as resolved by Lagrange until Weierstrass' 1858 and 1868 memoirs and Jordan's 1871 and 1872 notes would give two different ends to the network) ([32]). In the purpose of the questioning of its identity, we will be referring to this network as a "discussion" (between different authors as well as various theories) and we will designate it under the name of "the discussion on the equation to the secular inequalities in the planetary theory".

One of the main feature of this discussion is its origin which, as it has already been illustrated with Villarceau's 1870 note, was systematically traced back to the resolution given by Lagrange to problems of "small oscillations" and more precisely to the mechanical meanings he had associated to the algebraic nature of roots in stating a

[30] Given a family of forms $s\Phi-\Psi$, $|\Phi|\neq 0$, Jordan's theorem applied to the linear substitution $\Psi\Phi^{-1}$ shows that there exists a non singular substitution $U$ such that $U^{-1}(\Psi\Phi^{-1})U=J$ where $J$ is in canonical form. Thus $H(s\Phi-\Psi)K=sI-J$, where $H=U^{-1}$ and $K=\Phi^{-1}U$. The family of forms $s\Phi-\Psi$ can thus be reduced to its canonical form $sI-J$.

[31] The constitution of this network is therefore depending upon the choice of the 1874 controversy as a moment of reference for the bibliographic research. Alternative standpoints would lead to constitute different networks. For instance, taking as a point of departure the works of H. Poincaré on the three body problem discussed in A. Robadey's paper in this volume would lead to focus on astronomical references on the stability of the solar system (see also [Laskar 1992]). Taking as a moment of reference the works of Sylvester and E. Weyr on matrices in the 1880's would highlight some papers such as the 1850-1851 works in which Sylvester had introduced the terms matrices and minors and which play a minor role in our network as they were only quoted once by Darboux in 1874 (see [Brechenmacher 2006d]).

[32] From the standpoint of the 1930's modern algebra, the different problems appearing in the discussion would be considered as belonging to the theory of matrices and consisting in the reduction of a pair $(A,B)$ of matrices in $(D,I)$ where $A$ is symmetric, $B$ is definite symmetric, $D$ is diagonal and $I$ the identity matrix such as, for instance the linear differential systems with constant coefficients $BY'=AY$ which are related to the eigenvalue problem $AX=\lambda BX$. See [Gantmacher 1959, p. 311].

"general method" for deciding of a mechanical system's stability, and whose validity had remained unquestioned until Weierstrass proved in 1858 that the occurrence of "real, unequal and negative roots" was not a necessary and sufficient condition for the oscillations to remain bounded. When this method was extended in the 1770's to the description of the "secular inequalities" of the parameters determining the planetary orbits [Lagrange 1781, p. 125], the algebraic nature of roots was linked to the stability of the solar system and this new issue prompted Laplace's intervention in the discussion and his attempts to give a general demonstration of the bounded nature of the system. Laplace's 1789 demonstration was based on the "remarkable relations" *(1,2) = (2,1)* and *[1,2] = [2,1]* between the differential systems'coefficients ([33]):

$$0 = (I)\ \frac{d^2\xi}{dt^2} + (I,2)\ \frac{d^2\psi}{dt^2} + (I,3)\ \frac{d^2\varphi}{dt^2} + \&c + [I]\ \xi + [I,2]\ \psi + [I,3]\varphi + \&c.$$

$$0 = (2)\ \frac{d^2\psi}{dt^2} + (I,2)\ \frac{d^2\xi}{dt^2} + (2,3)\ \frac{d^2\varphi}{dt^2} + \&c + [2]\psi + [I,2]\xi + [2,3]\varphi + \&c$$

$$0 = (3)\ \frac{d^2\varphi}{dt^2} + (I,3)\ \frac{d^2\xi}{dt^2} + (2,3)\ \frac{d^2\psi}{dt^2} + \&c + [3]\varphi + [I,3]\xi + [2,3]\psi + \&c$$

As it is illustrated by some of the titles given to papers published within the discussion, such as Cauchy's 1829 "Sur l'équation à l'aide de laquelle on détermine les inégalités séculaires des planetes", Hermite's 1857 "Mémoire sur l'équation à l'aide de laquelle, etc. " or James Joseph Sylvester's 1883 "On the equation to the Secular Inequalities in the Planetary Theory" while the link of these papers with astronomy went no deeper than their titles ([34]), the mechanical works of Lagrange and Laplace had given a specific identity to an algebraic equation recognisable by the special nature of its roots and by the "remarkable relations" occurring in the linear systems related to this equation. It was thanks to this specific identity that Cauchy recognised in 1829 a formal analogy between various problems such as the small oscillations of mechanical systems, the rotation of a solid body or the classification of conics and quadrics [Cauchy 1829, p. 173] ([35]).

In a paper devoted to an extensive study of the discussion [Brechenmacher, 200?b], it has been shown that the special nature of the equation to the secular inequalities was closely related to a specific practice which resorted to the already mentioned polynomial quotients (*) $\dfrac{\frac{P_{1i}}{S}(x)}{x - s_j}$ in expressing the solutions of linear systems (such expressions were at first regarded as involving some equations obtained by elimination methods and, after some authors such as Cauchy or Jacobi made the determinants concept a basis for their analytic methods it gradually came to be regarded as involving successive sub-determinants $\frac{\Delta_{i+1}}{\Delta_i}$ extracted from *S*) ([36]). Not only did the "remarkable relations" property of the coefficients of systems we would nowadays designate as symmetric systems appeared as a consequence of the specific practice Lagrange had been developing while working out the initial conditions solutions to problems of small oscillations, a heritage of this practice handed down through the network of the discussion within different methods and various theoretical frameworks. Even though we would be tempted to detect in such a practice an origin of the method of transformation of a pair of symmetric matrices (*A,B*) into the pair (*D, I*) where *D* is a diagonal matrix and *I* the identity matrix, before the time of algebraic theories such as the 1930s theory of canonical matrix most of authors would not look upon this practice as a method within a theoretical framework. Related as it was to a specific equation, this practice had nevertheless taken on an algebraic identity within a consistent network and because the 1874 controversy highlighted some conflicting views on the nature of algebra (as when Kronecker was criticizing Jordan's algebraic organisation of the theory of forms) we shall look into this question more closely. Until the development of theoretical frameworks such as the ones Jordan and Kronecker were quarrelling about in 1874, it was above all a historical identity that characterised the discussion's algebraic nature. It was indeed in appealing to a corpus referring to earlier texts they themselves consequently contributed to develop that authors pointed out the specific nature of the equation to the secular inequalities. For instance, when Cauchy, in 1829, transcended the framework of analytical geometry he was originally interested in, it was thanks to the references he made to the mechanical works of Lagrange and Laplace in order to identify the algebraic practice he inserted in his own method for determining the principal axis of conics and quadrics with the result of generalising to *n* variables an analytical method originally devised for two or three variables.

[33] By modern standards, Laplace's proof, which was based on the conservation of the system's energy, was not valid because what he did was to use symmetry together with the differential equations to derive an equality that implied the solution had to be bounded as function of time. Then, from the form of the solutions to the differential equation –which were not correctly formulated for multiple roots – he inferred the reality of eigenvalues. See [Hawkins, 1975, p. 15].

[34] Cauchy' title was taken on at the last minute when Charles Sturm called his attention to the connection of his work in analytical geometry with secular perturbations; see [Hawkins 1975, p. 15].

[35] From the standpoint of modern algebra, Cauchy was interested in the transformation of a quadratic form in three variables into a sum of squares only. This problem also arose in the mathematical analysis of the rotational motion of a rigid body as studied by Lagrange in the 18th century. In Lagrange's work this problem was nevertheless not connected to the swinging string problems. For a description of Cauchy's work on the problem of the rotation of a solid body in connection to Lagrange's analytic reformulation of the solution given by Euler, see [Hawkins 1975, p.18].

From the outset of the text corpus to its two ends in Weierstrass' and Jordan's papers, it was an ambition of generality which was driving authors on making reference on a discussion they consequently joined themselves. It was with the aim of generalising to the "oscillations of an unspecified system of bodies" d'Alembert's investigations of a swinging a string loaded with two or three masses that Lagrange had been working out in 1766 the general *i.e.* polynomial practice at the origin of the discussion. Because of the generality he attributed to his description of a motion that Daniel Bernouilli had regarded as to be too irregular to be treated by analytic methods [Truesdell 1960, p. 156], Lagrange made the problem of small oscillations come out first among the examples of applications he gave to the "general principles" set in his 1788 *Mécanique analytique*. While generality was the main impetus given to the development of the discussion on the qualitative nature of characteristic roots, it nevertheless took on changing meanings between 1766 and 1874. To Lagrange's mind, the fact that his method would fail if multiple roots should occur did not restrict any of its generality because this method was resorting to implicit mechanical representations. Because it was known that the oscillations of a swinging string loaded with *n* masses could be mechanically represented as a combination of independent oscillations of *n* strings loaded with a single mass, differential systems had to be representable by combinations of independent equations $dy_i=kyi$ ($i=1, \ldots ,n$). In Lagrange's method algebraic roots could not be dissociated from their mechanical representations as periods of oscillations and the occurrence of multiple roots was therefore (wrongly) believed to be contradictory to the existence of *n* independent oscillations ; that is why they were considered to imply unbounded oscillations, an eventuality that had to be rejected "because of the nature of the problem" as it would be contradictory to the prerequisites that had been made on the bounded nature of the oscillations. D'Alembert, Lagrange and Laplace extended to the case of multiple roots the arguments, based upon physical considerations, they had developed for proving that the roots of the equation obtained by elimination method had to be real and negative : the problem concerned a swinging string and the displacements of the masses from the vertical must consequently remain small, but this would not be the case if the analytical solutions contain exponential $e^{\delta t}$ because then they would increase to infinity. The stakes in the implication of Lagrange's conclusion –the roots have to be real and unequal because the oscillations have to remain bounded – changed when the method was generalised to the secular inequalities in planetary theory. The stability of the solar system could not be taken for granted, Lagrange therefore pointed out that "it would be difficult, perhaps impossible, to determine the roots of the equation in general" [Lagrange 1766, p. 538] for it would mean demonstrating the reality and inequality of the roots of a very general $n^{th}$ degree polynomial while these roots could not be worked out in general as soon as *n* was be greater than five. Although he had worked out an effective computation for a fourth planets system and determined that the roots of the associated fourth degree equation were real, negative and unequal, Lagrange came to the conclusion that "one could wonder whether, by changing these values, perhaps equal or imaginary roots might occur. In order to eliminate all doubt, it would be necessary to demonstrate, in general, that the roots of the equation will always be real and unequal, whatever the values of the masses, provided only that they be positive. That is easy when the mutual action of only two planets is considered simultaneously, since the equation is only of the second degree, but this equation becomes more and more complicated and higher [in degree] as the number of planets increases" [Lagrange 1784, p. 316, translation T. Hawkins]. Laplace did not stay content with this numerical computation and his aim of devising a "fully general" demonstration that would not depend upon the approximate values assigned to the masses of the planets brought him to engage in the discussion. As has already been pointed out before, Cauchy's 1829 intervention was caused by his ambition to generalise a method he had devised for two or three variables in a geometric framework, the calculus of skew functions (determinants) he had appealed to for proving the principal axis theorem of conics and quadrics could not only be used to translate Lagrange's proof of the principal axis theorem for a rotational solid body but also the polynomial practice peculiar to the systems of *n* equations related to the equations to the secular inequalities ([37]). This generalization also induced some changes of perspectives on generality as some algebraic expressions emerging from the generalisation to $n$ variables of the method for changing rectangular coordinate systems could take a $\frac{0}{0}$ value in the case of multiple roots ([38]). The proof Cauchy had sketched for the reality of the roots was valid only if no successive equations $S_i(s)$ relating to the successive determinants extracted from $S(s)$ had a root in common. The occurrence of multiple roots now appeared as a singular case limiting the range of validity of an algebraic expression, it thus seemed necessary to introduce some particular methods for this singular case such as the infinitesimal argument Cauchy made use of in 1829 ([39]). But Cauchy did not rest content with this situation; aiming for a full homogeneous resolution as opposed to the singular cases which overburdened general polynomial methods, he developed his calculus of residues ([40]). The change of perspective on generality induced by this ideal of homogeneity would impulse the further developments of the discussion between 1830 and 1858, from Cauchy to Weierstrass and involving Jacobi, Borchardt, Hermite and Sylvester ([41]).

---

[36] As seen in note n°19, from the standpoint of modern algebra this practice could be considered as a method giving the general polynomial expressions of the eigenvectors of symmetric matrix *A* (such expressions are given by the columns of the matrix of cofactors computed from the polynomial matrix $A$-$\lambda I$). This formulation nevertheless induces some anachronisms not only because it resorts to modern notions or theories but also because it is implicitly related to complex methods of "transformations" (which includes some geometrical analogy such as the "symmetry" property of mechanical systems) which were extraneous to the practices used during the discussion. When he gave an integrable "form" to his systems, Lagrange did not resort to a method of "transformation" but to the computation of some mechanical parameters (the proper periods). The symmetric property of mechanical systems was a consequence of the practice devised by Lagrange for a direct computation of the solutions relating to the initial conditions from the characteristic equation, this practice, described in [Brechenmacher 200?b] could be considered as making use of dual orthogonlality from the standpoint of modern algebra.

[37] As Thomas Hawkins argued, "Cauchy succeeded in generalizing Lagrange's proof [..] by using the theory of determinants, a new mathematical tool he had brought to perfection in an earlier memoir [1815]; The reason

---

Cauchy probably decided to generalize the principal axis theorem of mechanics and quadratic surfaces was that it provided an occasion for him to apply his theory of determinants" [Hawkins 1977, p. 125].

[38] From the standpoint of modern algebra, the existence of the orthogonal transformation which diagonalized Cauchy's quadratic form depended upon the reality of the eigenvalues (*i.e.* the roots of the characteristic equation ) as well as upon the non existence of multiple roots.

[39] D'Alembert had already realized that his solutions to the swinging string problems might present some difficulties if the roots were not all distinct and he had appealed to some considerations on infinitesimals in order to handle the case of multiple roots. Such considerations would be retaken and developed later by Lagrange in 1766, Cauchy in 1829 and Sylvester in 1881.

[40] In 1870, Yvon-Villarceau would argue that Lagrange's mistake would originate in an abusive generalisation of the resolution of single equations of the $n^{th}$ order to systems of $n$ equations.As seen before, when Cauchy had composed his 1829's paper he was not particularly interested in systems of linear differential equations with constant coefficients. In the 1830's however, he became increasingly interested in the problem of deriving the properties of light from a theory of the small oscillations of a solid elastic medium. When he aimed at covering the case of multiple roots in 1839, he rejected this method and developed a uniform and completely general method of expression the solutions by using the calculus of residues he had introduced as soon as 1826 in order to cover the cases of multiple roots when integrating a single equation of the $n^{th}$ order with constant coefficients $\frac{d^n y}{dx^n} + a_1 \frac{d^{n-1} y}{dx^{n-1}} + ... + a_{n-1} \frac{dy}{dx} + a_n y = 0$. The algebraic resolution Cauchy was still relying to in 1825 resorted to an « analogy » between the polynomial expression $F(r)= r^n+a_1r^{n-1}+...+a_{n-1}r+a_n$ and a symbolical factorisation of the differential equation $(D-r_1)(D-r_2)...(D-r_n)y = - f(x)/a_0$, (the $r_i$'s being the roots of $F$). As $(D-r)y= f(x)$ implies $y = e^{rx} \int e^{-rx}f(x)dx$, Cauchy thus gave in 1825 the following expressions which is incorrect for multiple roots :

$$y = \frac{e^{r_n x}}{a_0} \int e^{(r_n - r_{n-1})x} (\int e^{(r_{n-1} - r_{n-2})x} (...\int e^{(r_2 - r_1)x} f(x)dx...)dx)dx$$

In 1826 he would give a solution whatever the multiplicity of roots thanks to the introduction of the calculus of residues : $y = \Re es \frac{\Phi(r)e^{rx}}{((F(r))}$ ($\Phi(r)$ being polynomial). See [Dahan Dalmedico 1992, p. 197].

[41] See [Hawkins 1977, p. 128-133] for some descriptions of the works of Jacobi and Borchardt and of the proof given by Dirichlet in 1846 (which would become an appendix to the third edition of Lagrange's *Mécanique Analytique* in 1853) to the fact that a state of equilibrium in a conservative mechanical system is stable if the potential function assumes a strict maximum value. Sylvester's works of 1850-1852 would lead him to introduce the notions of "matrix" and 'minors" (in which Darboux would see in 1874 an origin of Weerstrass' elementary divisors). See [Brechenmacher 2006d]. The methods developed by Sylvester and Cayley in the context of the

After he had found and stated the integrals [of mechanical systems], Lagrange came to the conclusion that, as the oscillations $x_1, \frac{dx_1}{dt}$ had to remain small if there were so at the origin, the equation could not possibly have multiple roots because in that case the integrals would increase to no limit over time. As he was dealing with the planetary secular variations in his *Mécanique céleste* Laplace repeated a similar affirmation. The same conclusion was mentioned by a number of author such as Poisson for instance. This conclusion is nevertheless groundless […] and the same result may be stated whether the roots of the equation f(s)=0 are distincts or not; if the homogeneity of the latter conclusion has not been reached before it is because this case [multiple roots] had always been dealt with by particular methods [Weierstrass 1858, p. 244, translation F.B.]([42]).

Even though ambitions of generality had been strongly linked to the algebraic identity of the discussion ever since its origin, it was actually the arithmetic nature of some processes of "transformations" of homogeneous forms of the second order" which supported Weierstrass' general and homogeneous conclusion. The only difference between the inertia law of the arithmetic theory of quadratic forms and Weierstrass' theorem was that the later was concerned with pairs of quadratic forms ($\Phi, \Psi$) and thus by the expressions $\Phi+s\Psi$ which, as it would be claimed later by Jordan, presented a polynomial and therefore algebraic nature. Whereas such terms as "forms" and "transformations" had been given an explicit mathematical definition in the arithmetic of quadratic forms in relation to the notion of equivalence relation that had been introduced by Gauss' 1801 *Disquitiones arithmeticae*, they pointed to various and mostly implicit meanings within the algebraic framework of the discussion. To Lagrange and Laplace minds, mechanical representations induced the existence of "integrable forms" for the differential systems of small oscillations and it was not by resorting to any idea of "transformation" that the independent equations were found but through the computation of the systems' mechanical parameters as well as the polynomial factorisations (*) of the characteristic equation ([43]). In Cauchy's 1829 paper, the "transformations of homogeneous functions" were related to geometrical meanings and the processes of changes of rectangular systems of coordinates for conics and quadrics. Even though it might seem natural nowadays to wonder about those matrices which, because of their multiple eigenvalues, might not be transformed into some diagonal forms, this question was actually irrelevant to the ways the terms "forms" and

---

development of the theory of invariants would be invested by Hermite in an arithmetical framework (quadratic forms, decomposition in four squares). On the birth of the theory of invariants. See [Parshall 1989 and 2006].

[42] Nachdem Lagrange die Form der Integral angegeben und gezeigt hat, wie die willkürlichen Constanten derselben durch die Anfangswerthe von $x_1, \frac{dx_1}{dt}$, u.s.w. stets unendlich klein bleiben, wenn sie es ursprünglich sind, auch die an, dass die genannte Gleichung keine gleiche Wurzeln haben dürfe, weil sonst in den Integralen Glieder vorkommen würden, die mit der Zeit beliebig gross werden könnten. Dieselbe Behauptung findet sich bei Laplace wiederholt, da wo er in der Mécanique céleste die Säcular-Störungen der Planeten behandelt, und ebenso, so viel mir bekannt ist, bei allen übrigen diesen Gegenstand behandelnden Autoren, wenn sie überhaupt den Fall der gleichen Wurzeln erwähnen, was z.B. bei Poisson nicht geschieht. Aber sie ist nicht begründet. [...], wenn nur die Function $\Psi$ stets negativ bleibt, und ihre Determinante nicht Null ist, was stattfinden kann, ohne dass die Wurzeln der Gleichung $f(s) = 0$ alle von einander verschieden sind ; wie man denn auch wirklich besondere Fälle der obige, Gleichungen, bei denen die Bedingung nicht erfüllt ist, mehrfach behandelt und doch keine Glieder von der angegebene Beschaffenheit gefunden hat.

[43] On the "mathematical interpretation of essential mechanical concepts" in Lagrange's analytical mechanic see [Panza 1992, p 205].

"transformations" were being considered at the times of Lagrange or Cauchy. Such kind of question was on the contrary a natural one from the standpoint of authors working in the 1850's on the arithmetic of quadratic forms –a framework which would be "generalised" to bilinear forms in the 1860's - as well as on the algebraic framework of Jordan's 1870 substitution theory ([44]). In the general resolution the latter gave in 1871 to the problem of the integration of linear differential systems with constant coefficients, such systems were algebraically "reduced" to a chain of "simple forms" related to the algebraic decomposition of the characteristic polynomial. At the two ends Jordan and Weierstrass had been giving to the discussion, to the different meanings and representations the term "form" had taken on within the methods of Lagrange, Laplace and Cauchy succeeded a mathematical theory whose subject was a "fully general" characterization of "forms". Should such a subject belong to arithmetic or algebra? Even though Jordan and Kronecker were referring in 1874 to the discussion corpus as a shared history relating to a specific practice they had in common and which consisted in investigating pairs of bilinear forms *(A, B)* by making use of the polynomial decomposition of $|A+sB|$, we shall see in the following section how some disciplinary ideals on algebra and arithmetic induced conflicting perspectives on generality.

### *5. Arithmetic generality vs algebraic generality.*

The purpose of the reorganisation Kronecker devised in 1874 for the theory of bilinear forms was to perform an arithmetical synthesis of different results that had been obtained in the 1860's ([45]). Although Kronecker had already been implicitly referring to the legacy of the works of Gauss and Hermite on the arithmetic of quadratic forms in 1866 – as when he preferred to make use of the term "form" to name what others would designate as a function ([Weierstrass, 1858]) or as a "polynom" ([Jordan, 1873]), his monthly communications to the Academy of Berlin during the winter of 1874 were aiming at an explicit generalisation of the arithmetic notion of "equivalence classes" from forms to networks of forms. "As an application of Arithmetic notions to Algebra", two bilinear forms or two networks of bilinear forms were designated as "equivalent" and as belonging to a same "class" when one could be linearly transformed into another ([46]). Some disciplinary ideals were coming along with this arithmetic orientation given to the theory and they expressed themselves in the criticisms Kronecker made of Jordan's statement that the identity of their canonical forms was a necessary and sufficient condition for two forms to be equivalent. According to Kronecker, despite being true, this proposition had to be rejected because it did not state any practical process for deciding of the equivalence ; it thus had to be distinguished from the "immediate possibility afforded by the theoretical criteria of equivalence to set a complete system of invariants" effectively computed from the form's coefficients as the result of the arithmetical process for computing the *g.c.d.*'s of the successive minors extracted from the polynomial determinant $|A+sB|$([47]).

---

[44] The origin of Christoffel's interest on "bilinear functions" (already mentioned in note n°9) was the generalisation he gave in 1864 of Weierstrass' 1858 theorem to the hermitian case for some purposes relating to Cauchy's theory light, Clebsch mechanical works of 1860, Hermite's results on biquadratic residues and in reference to Jacobi's 1857 memoir which gave a generalisation of the inertia law to "bilinear function" $\Sigma a_{ij}x_iy_j$. see [Mawhin 1981].

[45] According to Kronecker, this arithmetizing ambition had been stemming from discussions the latter had with *E. Kummer. On Kummer's ideal numbers see J. Boniface's paper in this volume.

[46] Two families of bilinear forms $s\Phi$-$\Psi$ and $s\Phi'$-$\Psi'$ are equivalent if one can be transformed into the other by (possibly different) non singular linear transformations of the $x$ and $y$ variables (where $\Phi = \sum_{i,j=1}^{n} A_{ij}x_ix_j$ and $\Psi = \sum_{i,j=1}^{n} B_{ij}x_ix_j$ ).

> In the arithmetical theory of forms, we actually have to content ourselves with indicating a method in order to decide of the equivalence […](cf. Gauss : Disquitiones arithmeticae, Sectio V [...]). Although this method might induce to resort to reduce forms, we should make clear that in the arithmetic theory of forms such reduced forms would be endowed a relevance completely different from the one they would assume in algebra. Because of their own nature, the invariants of equivalent forms must be obtained as the results of arithmetical processes performed on the forms' coefficients, and it is not surprising that such processes, despite being directly defined, can not be represented explicitly as the results of some arithmetical operations; as it so happens for most of the notions of arithmetic such as the simple notion of greater common divisor [Kronecker 1874b, p. 415, translation F.B.] ([48]).

The ideal of effectivity, which the historiography had been usually linking to the 1880's arithmetic theory of algebraic magnitudes ([49]), had thus already been strongly expressed on the occasion of the 1874 controversy when Kronecker blamed "litteral expressions" such as Jordan's canonical form for needing an algebraic decomposition of the characteristic determinant and for which no practical process could therefore be given "in general" as soon as the polynomial degree would exceed five.

Throughout the 1874 controversy, Jordan was retorting to Kronecker's assaults by claiming the greater generality and simplicity of his method. Way off the naïve simplism caricatured by Kronecker, Jordan's ideal of simplicity was linked to a practice of "reduction" of "general problems" into chains of sub problems. It supported a criticism of Kronecker's 1868 characterization of singular pairs of bilinear forms as having failed to find the "true

---

[47] Consider a family of bilinear forms $s\Phi\text{-}\Psi$ and let $S(s)$ denote the determinant $|s\Phi\text{-}\Psi|$, the greatest common divisor of all the first minors of $S(s)$ considered as polynomials in $s$ is denoted by $S_1(s)$. Similarly, $S_2(s)$ is defined as the greatest common divisor of all the second minors of $S(s)$ and so on. Then $S_i(s)$ divides $S_{i-1}(s)$ and if $E_i(s)$ denotes the polynomial $S_{i-1}(s)/S_i(s)$ then $E_i(s)$ divides $E_{i-1}(s)$. Thus, $S(s)$ differs from the product of the $E_i(s)$ by a constant and if $s_1, s_2, \ldots, s_k$ are the distinct roots of $S(s)$ then

$$E_i(s) = c_i \prod_{j=1}^{n} (s - s_j)^{m_{ij}}$$

where $c_i$ is constat and the $m_{ij}$ are positive integers or zero. Each factor $e_{ij} = (s - s_j)^{m_{ij}}$ with $m_{ij}>0$ is an elementary divisor of $S(s)$.

[48] In der arithmetischen Theorie der Formen muss man sich freilich mit der Angabe eines Verfahrens zur Entscheidung der Frage der Aequivalenz begnügen und das betreffende Problem wird deshalb auch ausdrücklich in dieser Weise formuliert (cf. Gauss : Disquitiones arithmeticae, Sectio V [...]) Das Verfahren selbst beruht auch dort auf dem Uebergange zu reductiren Formen : doch ist dabei nicht zu übersehen, dass denselben in den arithmetischen Theorien eine ganz andere Bedeutung zukommt als in der Algebra. Da nämlich die Invarianten äquivalenter Formen dort ihrer Natur nach nur zahlentheoretische Functionen der Coëfficienten sind, so kann es nicht befremden, wenn dieselben zwar direct definiert aber nicht explicite sondern nur als Endresultate arithmetischer Operationen dargestellt werden können ; denn ganz ähnlich verhält es sich mit den meisten arithmetischer Begriffen, z.B. schon mit jenem einfachsten Begriffe des grössten gemeinsamen Theilers.

[49] The arithmetical aspects of polynomials already played an important role in Kronecker's 1850-1870 work on the solvability of equations and would later be essential in his 1882 arithmetic theory of algebraic magnitudes and his conception of a "Rationalsbereich" which was based on polynomial forms as an alternative to Dedekind's fields.

reduced forms" which had to be simplest links of the chain of reductions with no possibility of further simplification [Jordan 1874b, p. 13]. Jordan's practice of reduction originated in the methods the author had devised in his researches on substitutions groups in the 1860's. The purpose of these researches to whose the *Traité des substitutions et des equations algébriques* gave a synthesis in 1870 was the general investigation of the various types of equations that could be solved by radicals. In order to handle this very general problem, the *Traité* developed a "machinerie", an "enormous récurrence on the degree $n$ of the equation" as Jean Dieudonné would describe it [Dieudonné 1970, p. 168], which reduced the "gender" of a substitution group from the general to the particular. The investigation of general soluble groups was therefore being reduced to the analysis of successive particular groups such as "transitive", "primitive", "linear" or "symplectic" groups corresponding to the simplest links of Jordan's chain of reduction. Among others, the linear group and its properties such as the theorem stating the reduction of linear substitutions to their "simplest" or "canonical" forms, "originated" in the practice of reduction Jordan made use of in his investigations of the 1860s and it was only afterwards that a theoretical organisation would be given to such properties in the *Traité* of 1870 ([50]).
In 1871, when Jordan came to answer to the questions Yvon-Villarceau had been asking to the geometers of the Academy, he did resort on his practice of reduction in bringing general systems of linear equations down to a sequence of "simplest forms" associated to the decomposition of the characteristic equation into its simplest (*i.e.* linear) factors. Requiring the resolution of a general algebraic equation, Jordan's canonical reduction showed an abstract nature as it did not actually provide any practical resolution of the problem. It was thanks to the practices –such as canonical reduction – he had originally devised as methods for group theory that Jordan succeeded in extending the range of his investigations to subjects such as differential equations (1871-1878), the theory of forms (1872-1875) as well as arithmetic and number theory (1878-1908) ([51]). These practices did not come alone in the applications and Kronecker's criticisms highlighted some of the ideals - such as simplicity and abstraction - which were walking along with them.
As an outcome of the 1874 controversy, in the two successive memoirs he would publish in 1878 and 1879, Frobenius would give to the theory of bilinear forms an organisation which would remain unchanged until the 1920s. Although Frobenius would elaborate his theory upon an orientation resuming the ambition Jordan had in 1874 to develop a single theory for both bilinear forms and linear substitutions, he would take up arithmetizing tendency of the theoretical organisation devised by Kronecker in focusing on rational invariant computations ([52]).

[50] From a contemporary standpoint, in order to characterize which primitive equations were solvable by radicals, Galois had asserted that the degrees of such equations were of he form $p^n$, $p$ prime and that the corresponding group $g$ of permutations had to be a solvable subgroup of the linear group. As opposed to his predecessors who often concentrated upon projective linear substitutions, Jordan made the consideration of homogeneous linear substitutions fundamental in his 1867 investigations on the determination of all the irreducible equations of a given degree which were solvable by radicals. It was for the purpose of establishing the three general types of solvable subgroups of the group of linear substitutions in two variables that he stated in 1868 the three kinds of "canonical forms" such a linear substitution $S$ could be reduced to depending of the nature of the roots of $|S-kI|\equiv 0|p]$. In 1870 the general canonical form theorem for $n$ variables played a major role in Jordan's method of building up solvable subgroups from their composition series which involved determining linear substitutions which commuted with a given substitution $S$. For a detailed analysis of the role played by canonical forms in Jordan's investigations on solvable groups as well as on the evolution of the role played by linear groups between 1870 and 1900 see [Brechenmacher 2006a].
[51] See [Brechenmacher 200?a] on the relation between this practice of reduction and the Jordan-Hölder theorem.

> By means of rational operations one can therefore determine whether given forms $A=sA_1+A_2$, $B=rB_1+B_2$ are equivalent or not. By contrast, all proofs known to me for the theorem of Herr Weierstrass involve irrational operations, for they are based on the transformation of $A$ into reduced forms whose coefficients depend upon the roots of the equation $|A|=0$. Hence I long ago had proposed to myself the problem of finding a proof for that theorem in which only rational operations occur. [Frobenius 1879, p. 483, translation T. Hawkins 1977, p.149].

Frobenius' success in giving a rational proof to Weierstrass' theorem in using the invariant factors introduced by Kronecker was due to his generalisation of Smith's lemma for ordinary integers [Smith 1861] to forms $A$-$sB$ with coefficients which were polynomials in $s$ and could be treated by analogy by arithmetical methods such as the Euclidean algorithm ([53]). The theory would thus revolve around the elementary divisor theorem for the

---

[52] As it has been discussed earlier, Jordan pointed out in 1873 that the transformation of bilinear forms could be seen as the composition of linear substitutions. Frobenius 1878-1879 papers developed a symbolical calculus on bilinear forms and added to the operations of addition and scalar multiplication the multiplication $AB=\sum\frac{\partial A}{\partial y_\gamma}\frac{\partial B}{\partial x_\gamma}=\sum c_{\alpha\beta}x_\alpha y_\beta$ where $c_{\alpha\beta}=\sum a_{\alpha\gamma}b_{\gamma\beta}$. The transformation of the form $A$ by means of linear substitutions $x_\alpha=\sum_\beta p_{\alpha\beta}X_\beta$, $y_\alpha=\sum_\beta q_{\alpha\beta}Y_\beta$ thus corresponded to the product $P'AQ$ , $P$ and $Q$ being bilinear forms identified with the linear substitutions ($P'$ denoting the transpose of $P$) :

$$P=\sum p_{\alpha\beta}x_\alpha y_\beta\text{, }Q=\sum q_{\alpha\beta}x_\alpha y_\beta$$

Frobenius thus gave some symbolical expressions to the equivalence relations introduced by Kronecker in 1874. This symbolical calculus reflected an heritage of some English authors such as Cayley and Sylvester whose works had been called to Frobenius attention by Darboux's 1874 paper and the same preoccupations both authors had in the 1870's for the Pfaffs' probem, see [Brechenmacher 2006a].

[53] From the standpoint of modern algebra, Kronecker's invariant factors can be introduced for any Euclidean ring while Frobenius' derivation of the elementary divisors is valid for any field. Applied to a bilinear form $A$ with integral coefficients, the lemma Smith had stated in 1861 asserts that $A$ is equivalent to (that is $A$ can be transformed into the following form by unimodular substitutions whoses determinants are +- 1)

$$F=f_1x_1y_1+f_1f_2x_2y_2+....+f_1f_2...f_lx_ly_l$$

where $l$ is the rank of $A$. Let $e_\lambda=f_1f_2..f_\lambda$ and $d_i=e_1e_2...e_i$ then $d_i$ and $e_i=d_i/d_{i-1}$ are the analogues of the polynomial invariants $S_i(s)$ and $E_i(s)$ of Kronecker's theory (see note n°48).

The derivation of Weierstrass' theorem from Smith's lemma was howewer not immediate as the latter implied that if a form $A$ had coefficients that were polynomials in $s$ then $PAQ=F$ where $P$ and $Q$ have coefficients which were polynomials in $s$. It was therefore necessary to prove that $P$ and $Q$ could be replaced by forms $P_0$ and $Q_0$ with scalar coefficients when $A$ was polynomial of first degree in $s$ (such as for the forms' similarity expressed by $A$-$sI$), Thomas Hawkins has pointed out that Frobenius' proof of the existence of $P_0$, $Q_o$ illustrated the operationally of his symbolical calculus of form which would later constitute an important method in matrix theory.

For some perspectives about the tendency to view objects arithmetically that were not originally regarded as within the province of arithmetic -such as Euclidean algorithm or uniqueness of factorisation of polynomials - reflected in this generalisation as inspired by Gauss'1832 introduction of "Gaussian integers" to establish the law of biquadratic reciprocity and Dededind's 1857 theory of higher congruences see [Hawkins 1977, p.150]. The objects of Dedekind study were polynomials with integral coefficients taken modulo a prime $p$ which adhered to equivalence of polynomial forms from which Jordan's canonical form would follow readily as a corollary thereby losing its theorem status ([54]). Following Weierstrass' 1868 paper, Frobenius would consider canonical or normal forms only for the purpose of proving the existence of a form $A+rA'$ of first degree whose coefficients belonged to some field and which possessed some prescribed elementary divisors $(\Phi(r))^{\varepsilon}, (\Phi_1(r))^{\varepsilon'}, \ldots$ Frobenius constructed such a form as a sum of forms, each with its own set of variables, $\varepsilon$ of them corresponding to $\Phi(r)$, $\varepsilon'$ to $\Phi_1(r)$ etc. Corresponding to $\Phi_i(r)$ where forms $r(x_1y_1+\ldots+x_\alpha y_\alpha)+y_1(a_1x_1+\ldots.+a_\alpha x_\alpha)-(x_1x_2+x_2y_3+\ldots+x_{\alpha-1}y_\alpha)$ with determinant

$$\begin{vmatrix} r+a_1 & -1 & 0 & \ldots & 0 \\ a_2 & r & -1 & \ldots & 0 \\ a_3 & 0 & r & \ldots & \ldots \\ \ldots & \ldots & \ldots & \ldots & 0 \\ a_\alpha & 0 & 0 & \ldots & r \end{vmatrix}$$

Frobenius concluded his 1879 memoir by explaining how the "canonical form that Jordan had stated in his *Traité des substitutions*" was thus obtained in the particular "case where he elements $a_{\alpha\beta}$ of the system $A$ are polynomials in $r$ with integral coefficients and two such functions are not considered to be distinct when their corresponding coefficients are congruent relative to a prime number modulus $p$ […] and if use is made of the complex numbers introduced by Galois" [Frobenius 1879, p…., translation T. Hawkins 1977, p. 153]. Even though the successive works of the Berliners on bilinear forms had developed some different methods - while Weierstrass would appeal to the successive algebraic decompositions of a polynomial determinant and its minors, Kronecker would compute *g.c.d.s* and arithmetical decompositions whereas Frobenius' theory would be based on a symbolical calculus on forms - these different methods would resort to a same practice. This practice, which consisted in resorting to polynomial invariants in order to handle mathematical "forms", was not restricted to a single method, theory or discipline and was strongly linked to some cultural aspects peculiar to the Berliners such as some modalities of handling generality in mathematics. As it would not appeal to any process of "transformation" but would resort to some polynomial computations, this practice was also inducing a specific static way of thinking about "forms", it was actually not until the computation had been done that the set of invariants that had thus been worked out would eventually be represented by some normal forms which would be nothing more than a specific way of representing a determinant.

On the opposite of the static nature of invariant computations, the practice Jordan had developed for the reduction of a general problem to a chain of simpler problem was inducing some dynamic ways of thinking about "transformations" and "reductions". The operationality of the representation Jordan had developed for linear substitutions allowed to "see" not only how the "indices" could be put into some sub groups in relation to the the analogy with elements of the theory the theory of number such as the Euclidean algorithm from which the basic arithmetical properties were then derived.

[54] Until 1907 when de Séguier presented in note to the *Comptes Rendus* Jordan's canonical form as the basis of a new organisation for the "theory of matrices", Jordan's result would only be considered as a theorem within the framework of group theory, on the posterity of Jordan see [Brechenmacher 2006a].

successive polynomial decomposition of the determinant but also the actions of the substitution on these sub groups ([55]) :

> The distinct functions $y_0, y'_0,...;...;y_p, y'_p,...;...$ can be chosen as independent indices in place of an equal number of primitive indices $x, x', ...x^{n-1}$. Let $m$ be the number of such functions, the substitution would thus be reduced to the form
>
> $$A=\left|\begin{array}{ll} y_0 & K_0 y_0 \\ y'_0 & K_0 y'_0 \\ \dots & \dots \\ y_1 & K_1 y_1 \\ \dots & \dots \\ x^m & a_1^m x^m + b_1^m x^{m+1} + \dots + c_1^m x^{n-1} + d_1^m y_0 + e_1^m y'_0 + \dots + f_1^m y_1 + \dots \\ x^{m+1} & a_1^{m+1} x^m + b_1^{m+1} x^{m+1} + \dots + c_1^{m+1} x^{n-1} + d_1^{m+1} y_0 + e_1^{m+1} y'_0 + \dots + f_1^{m+1} y_1 + \dots \\ \dots. & \dots \\ x^{n-1} & a_1^{n-1} x^m + b_1^{n-1} x^{m+1} + \dots + c_1^{n-1} x^{n-1} + d_1^{n-1} y_0 + e_1^{n-1} y'_0 + \dots + f_1^{n-1} y_1 + \dots \end{array}\right|$$
>
> [Jordan 1870, p. 117, translation F.B.].

In representing simultaneously the indices and the actions of the substitution on these indices, Jordan's notation was depicting the successive steps of the practice of reduction itself. The form in the above quotation contained the substitution $C$ on which the reduction had to be repeated

$$c=\left|\begin{array}{ll} x^m & a_1^m x^m + b_1^m x^{m+1} + \dots + c_1^m x^{n-1} \\ x^{m+1} & a_1^{m+1} x^m + b_1^{m+1} x^{m+1} + \dots + c_1^{m+1} x^{n-1} \\ \dots & \dots \\ x^{n-1} & a_1^{n-1} x^m + b_1^{n-1} x^{m+1} + \dots + c_1^{n-1} x^{n-1} \end{array}\right|$$

This process of reduction had to be repeated until all the indices would be "grouped" in "distinct series $Y_0, Z_0, u_0, ... ; Y'_0, Z'_0, ... ; ...$ " on which $A$ would act by the "simple law" $K_0Y_0, K_0(Z_0+Y_0), K_0(u_0+Z_0), ...; ...$ associated to the canonical form :

$$\left|\begin{array}{ll} y_0, z_0, u_0, \dots, y'_0, \dots & K_0 y_0, K_0(z_0+y_0), K_0(u_0+z_0), \dots, K_0 y'_0 \\ y_1, z_1, u_1, \dots, y'_1, \dots & K_1 y_1, K_1(z_1+y_{10}), K_1(u_1+z_1), \dots, K_1 y'_1 \\ \dots\dots\dots & \dots\dots\dots \\ \upsilon_0 & K'_0 \upsilon_0, \dots \\ \dots\dots\dots & \dots\dots\dots \end{array}\right|$$

### *Conclusion.*

We shall now come to some conclusions that may be drawn from the conflicting perspectives on generality relating to the two practices opposed by Jordan and Kronecker in 1874. In a word, while Jordan, on the one hand, criticized the lack of generality of Kronecker's invariant computations because they would not state explicitly the simplest reduction of pairs of forms, Kronecker considered Jordan' canonical form as a "formal notion" with no "objective meaning" thereby displaying a "generic" nature and failing to reach a true generality. It was in the first place the "general" solution they achieved for different problems that had been handled in the past by authors such as Lagrange, Laplace, Cauchy or Hermite, and furthermore their capacity to consider these various problems as a single "general" problem of "transformation" of pair of "forms" that prompted some connections between Jordan's practice of canonical reduction and the practice of invariants computations used in the development of the theory of bilinear forms in Berlin. The reference to a common history therefore played key role in the controversy. Not only did Kronecker, with the aim of criticizing Jordan's canonical form, stress a history of what the historian T. Hawkins would later refer to as the "generic reasoning" in the $18^{th}$-$19^{th}$ centuries algebra, but the reference to the "discussion on the equation to the secular inequalities in the planetary theory" played also a major role in identifying a specific practice that consisted in expressing the solutions of linear equations as polynomial factors of their characteristic equation. Even though they were sharing this traditional polynomial practice of generality, Jordan and Kronecker melted it with practices of their own in order to design the two methods related to their conflicting views on the theoretical organisation of a "general theory of forms". On the one hand, in the algebraic organisation Jordan gave to the theory, transformations resulted from the action of some linear substitutions groups and, in order to achieve "general results" on forms, underlying substitutions had to be reduced to their "simplest canonical forms" depending on the nature of the linear group the substitutions were belonging to. On the other hand Kronecker blamed algebraic methods for their tendency to think in term of the "general" case with little attention given to difficulties that might be caused by assigning specific values to the symbols whereas Weierstrass' invariants were considered as a model of "truly general" development. Even though Jordan's canonical forms could not be charged with an indictment of "so called generality" as they were actually giving a solution whatever the multiplicity of the characteristic roots, it was the algebraic nature of the practice of reduction they were related to that prevented them from reaching a general and theoretical level because, Kronecker argued, one shall not mistake algebraic methods for the "general notions" relating to arithmetic it was algebra duty to serve. Kronecker appealed to the tradition of Gauss on behalf of his claim that the theory of forms should be considered as belonging to arithmetic and should consequently focus on the characterisation of equivalence classes in establishing arithmetical invariants thanks to some effective procedures such as *g.c.d.*s computations. As long as they could not be effectively computed because they resorted to the solution of "general" algebraic equations, explicit resolutions such as Jordan's canonical form had thus to be rejected because of their formal nature.

[55] The $K_i$ designating the characteristic roots. From the standpoint of modern algebra, a vector space is decomposed under the action of $A$ as a sum of a stable subspace and his complementary space. On Jordan's demonstration see [Brechenmacher, 2006a, p.167-187].

On the opposite of Frobenius' 1879 theoretical organisation, the "Jordan canonical form theorem" would be referred to as a central result of the "theory of matrices" in most of the many treaties that would be published in the 1930's. This theorem would actually relate to two canonical forms ([56]):

[56] The two matrices given in this example both relate to the minimal polynomial $\lambda^8+\alpha_1\lambda^7+\ldots+\alpha_7\lambda+\alpha_8=(\lambda-\lambda 1)^2(\lambda-\lambda_2)^3(\lambda-\lambda_3)(\lambda-\lambda_4)$.

$$A=\left\|\begin{matrix} \lambda_1 & 1 & 0 & & & & & 0 \\ 0 & \lambda_1 & 0 & & & & & \\ & & \lambda_2 & 1 & 0 & & & \\ & & 0 & \lambda_2 & 1 & & & \\ & & 0 & & \lambda_2 & 0 & 0 & 0 \\ & & & & 0 & \lambda_3 & 0 & 0 \\ & & & & & & \lambda_4 & 1 \\ & & & & & & & \lambda_4 \end{matrix}\right\|,\ B=\left\|\begin{matrix} 0 & 0 & \dots & & & \dots & 0 & -\alpha_8 \\ 1 & 0 & \dots & & & & \dots & -\alpha_7 \\ 0 & 1 & 0 & \dots & \dots & & & -\alpha_6 \\ \dots & 0 & 1 & 0 & 0 & \dots & & \\ & \dots & 0 & 1 & 0 & 0 & \dots & \dots \\ & & \dots & 0 & 1 & 0 & 0 & \dots \\ & & & & 0 & 1 & 0 & -\alpha_2 \\ 0 & 0 & \dots & & \dots & 0 & 1 & -\alpha_1 \end{matrix}\right\|$$

On the one hand the Jordan canonical matrix $A$ would be considered as the simplest form in connection to the maximal decomposition of a matrix, on the other hand the "rational canonical form" $B$ would be obtained as the result of effective procedures. While laying an emphasis on canonical forms the theory of canonical matrices would put to the foreground a method of decomposition resorting to an ideal of simplicity closed to the one Jordan had been appealing to in 1874:

> The theory of canonical matrices is concerned with the systematic investigation of types of transformation which reduce matrices to the simplest and most convenient shape. [Aitken and Turnbul, 1932, 1].

As it is illustrated by the examples given below from the *Theory of canonical matrices* written by C.C. Turnbull and A.G. Aitken, this method would resort to some operatory processes relating to the representation given to matrices in order to decompose matrices into partitions of submatrices :

> The theory of canonical matrices is concerned with the systematic investigation of types of transformation which reduce matrices to the simplest and most convenient shape […]. It is convenient to extend the use of the fundamental laws of combination for matrices to the case where a matrix is regarded as constructed not so much from elements as from submatrices, or minor matrices, of elements. For example, the matrix
>
> $$A=\left[\begin{array}{cc:c} 1 & 2 & 3 \\ 4 & 5 & 6 \\ \hdashline 7 & 8 & 9 \end{array}\right]$$
>
> can be written
>
> $$A=\begin{bmatrix} P & Q \\ R & S \end{bmatrix}.$$
>
> 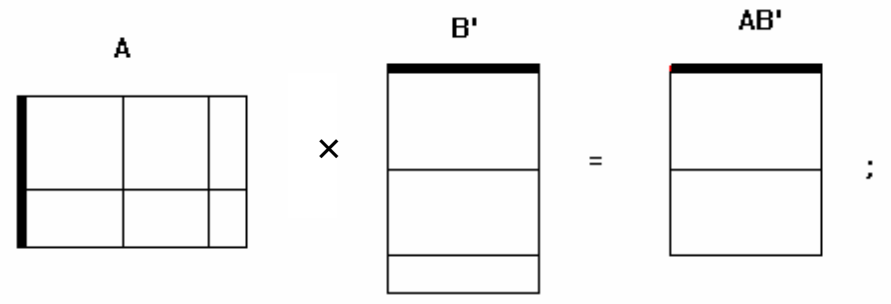
> 
>
> […]
>
> It is a help to form a staircase graphs of these chains as follows:

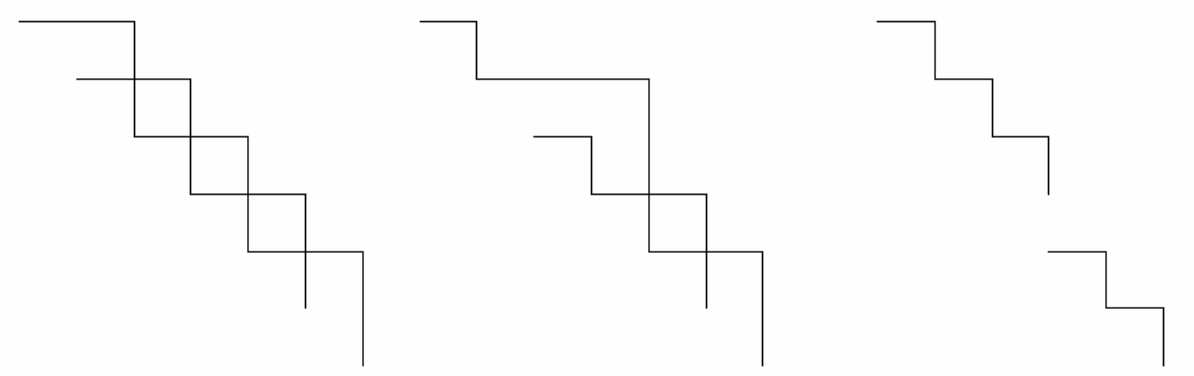

Fig. I Fig. 2 Fig.3

$$\begin{bmatrix} \alpha & 1 & . & . & & & \\ & \alpha & 1 & . & & & \\ & & \alpha & 1 & & & \\ & & & \alpha & & & \\ & & & & \alpha & 1 & . \\ & & & & & \alpha & 1 \\ & & & & & & \alpha \end{bmatrix}.$$

[Turnbul et Aitken, 1932, 1-20].

Between 1874 and 1930, the tension between canonical forms and invariants which has been discussed in this paper, would play a major role in the complex history of the practices that would later come to give an operatory dimension to the matrical representation as the basis of a method melting the polynomial generality of some invariants to some new processes such as an algebraic symbolical calculus, a combinatory of submatrices, an arithmetic of rows and columns as well as the geometric decomposition of a vector space in stable subspaces, thereby depicting how a general problem could be reduced into a chain of simpler problems.

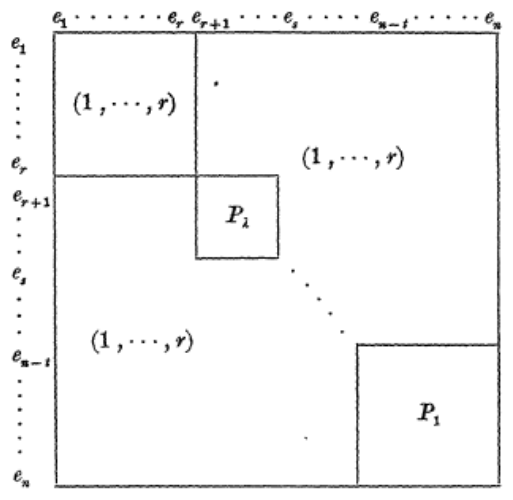

A careful study of the dynamic of this tension between canonical forms and invariants practices was used as a preliminary to a wider historical understanding of the history of linear algebra in our doctoral thesis. On the fringe of a dominant theory of bilinear forms which would resort to practices of invariants computations and which, as it has been shown by Thomas Hawkins would globalize some ideals on generality developed in Berlin such as the opposition of "true generality" and generic reasoning, some authors such as Jordan, Henri Poincaré, Eduard Weyr, Theodor Molien, Kurt Hensel, William Burnside, Leonard Dickson or Léon Autonne would handle general problems in developing some practices of reductions to canonical forms which would resort to some operatory processes on imagery representations. Some of these practices would not resort to what we would identity today as the framework of linear algebra – such as the calculus of tables developed mainly in an arithmetical framework by authors such as Hermite, Jordan, Poincaré and Chatelet – and studying how they were fitting in some networks of references not only highlights how these practices were related to some modalities of handling generality peculiar to these networks but it also raises some issues about disciplines and communities formations, evolutions and connections in relation to the ambitions of "generalisations" which are of particular importance for the history of linear algebra as it has been portrayed in this paper for the case of the network of the discussion on the equation to the secular inequalities of planetary theory.

## Annex.
## The discussion on the equation to the secular inequalities on planetary theory.

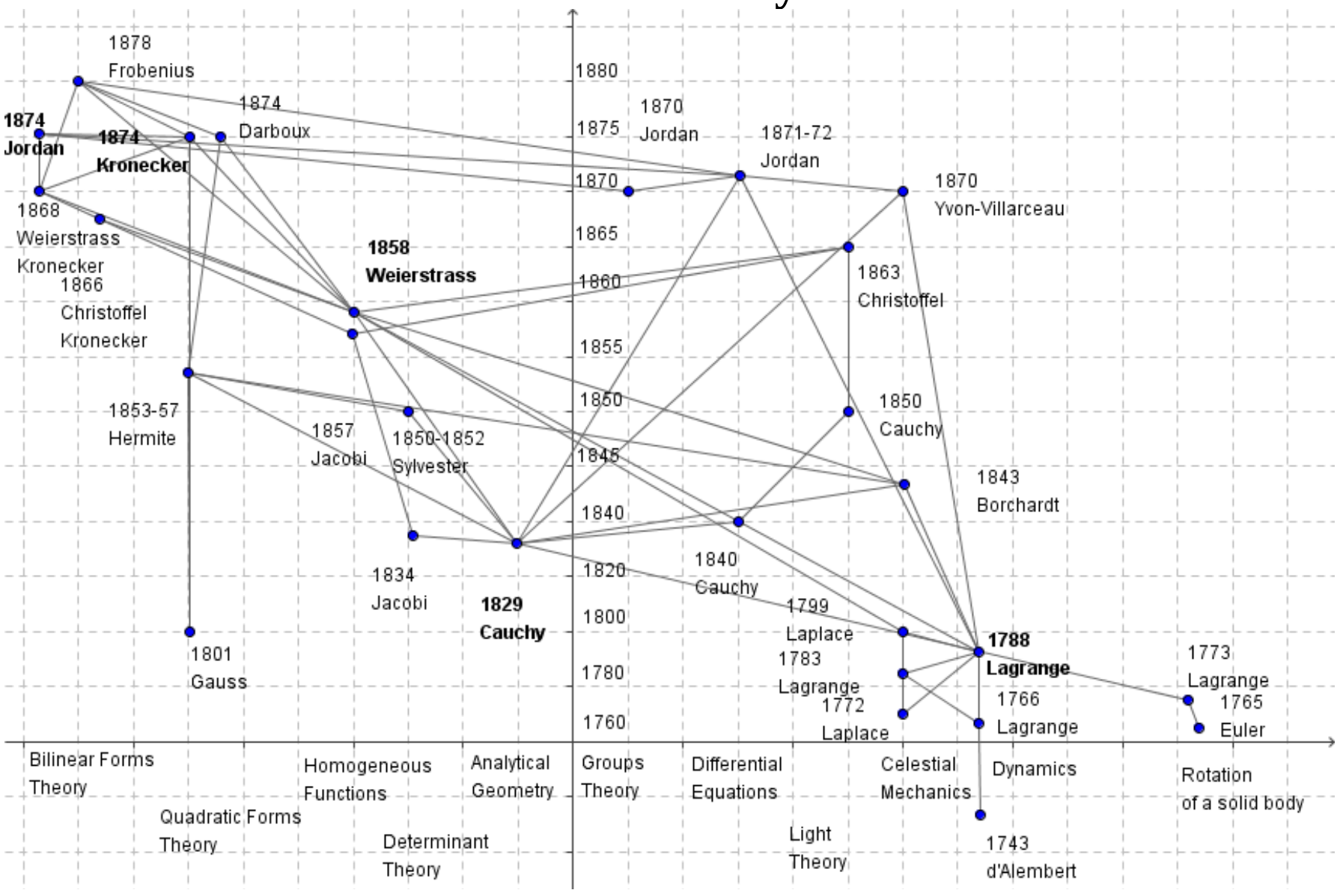

## Bibliography.